\documentclass[journal,twoside,web,onecolumn]{ieeecolor}
\usepackage{generic}
\usepackage{cite}
\usepackage{amsmath,amssymb,amsfonts}
\usepackage{graphicx}
\usepackage{textcomp}

\usepackage{graphicx}%
\usepackage{multirow}%

\usepackage{mathrsfs}%

\usepackage{xcolor}%
\usepackage{textcomp}%
\usepackage{manyfoot}%
\usepackage{booktabs}%

\usepackage{algpseudocode}%

\newtheorem{definition}{Definition}
\newtheorem{theorem}{Theorem}[section]
\newtheorem{proposition}[theorem]{Proposition}
\newtheorem{lemma}{Lemma}

\allowdisplaybreaks[4]

\begin{document}
\title{Open-Loop and Closed-Loop Strategies for Linear Quadratic Mean Field Games: The Direct Approach}
\author{Yong Liang, Bing-Chang Wang, and Huanshui Zhang
\thanks{ Corresponding author Bing-Chang Wang.  This work was supported by the National Natural Science Foundation of China under Grants 61922051, 62122043, 62192753, 62250056,  Natural Science Foundation of Shandong Province for Distinguished Young Scholars (ZR2022JQ31), 
	Major Basic Research of Natural Science Foundation of Shandong Province (ZR2021ZD14, ZR2020ZD24). }
\thanks{ Yong Liang is with the School of Information Science and Engineering, Shandong Normal University, Jinan, 250358, Shandong, P. R. China (e-mail: yongliang@sdnu.edu.cn). }
\thanks{Bing-Chang Wang is with the School of Control Science and Engineering, Shandong University, Jinan, 250061, Shandong, P. R. China (bcwang@sdu.edu.cn).}
\thanks{Huanshui Zhang is with 
the College of Electrical Engineering and Automation, Shandong University of Science and Technology, Qingdao, 266590, Shandong, P. R. China (e-mail: hszhang@sdu.edu.cn).}}

\maketitle

\begin{abstract}
This paper delves into studying the differences and connections between open-loop and closed-loop strategies for the linear quadratic (LQ) mean field games (MFGs) by the direct approach. The investigation begins with the finite-population system for solving the solvability of open-loop and closed-loop systems within a unified framework under the global information pattern. By a comprehensive analysis through variational methods, the necessary and sufficient conditions are obtained for the existence of centralized open-loop and closed-loop Nash equilibria,  which are characterized by the solvability of a system of forward-backward stochastic differential equations and a system of Riccati equations, respectively. The connections and disparities between centralized open-loop and closed-loop Nash equilibria are analyzed. Then, the decentralized control is designed by studying the asymptotic solvability for both open-loop and closed-loop systems. Asymptotically decentralized Nash equilibria are obtained by considering the centralized open-loop and closed-loop Nash equilibria in the infinite-population system, which requires a standard and an asymmetric Riccati equations. The results demonstrate that divergences between the centralized open-loop and closed-loop Nash equilibria in the finite-population system, but the corresponding asymptotically decentralized Nash equilibria in the infinite-population system are consistent. Therefore, the choice of open-loop and closed-loop strategies does not play an essential role in the design of decentralized control for LQ MFGs.
\end{abstract}


\section{Introduction}

Over the past decade, mean field games (MFGs) have captivated the collective interests of the mathematical and control engineering communities, and a series of research results have been widely achieved  across a broad spectrum including system control, applied mathematics and economics \cite{b102301,b102302,b102303}. In large-population systems, it is impractical to implement centralized control strategies due to the limitations of sensing, computing  and communication capabilities. However, the MFG theory provides an efficient theoretical scheme for handling such complexity. The most salient feature of mean field models is the interactive weakly-coupling structure across considerable players: the individual influence on the entire system is negligible, but the aggregate impact of overall populations to a given player is significant and cannot be ignored. By mean field approximations, the decentralized $\varepsilon$-Nash equilibrium can be designed only using its own individual state and some quantity computed offline accordingly, where $\varepsilon\rightarrow0$ as $N\rightarrow \infty$ with $N$ is the number of players in the large-population system. Therefore, the computational complexity can be  significantly reduced. The past few years  have showcased a plethora of successful MFG applications in various fields, such as dynamic production planning \cite{b102402}, vaccination games \cite{b102401}, opinion dynamics \cite{b102403,b013112}, dynamic collective choice \cite{b013110,b013111}, etc.     

MFGs germinated independently from the seminal works of Lasry and Lions \cite{b102404}, as well as Huang et al. \cite{b102405}. Readers are referred to the surveys \cite{b102406,b102407} for a useful overview. According to the timing of using mean field approximations, the development of MFG theory has largely bifurcated into two routes: the fixed-point approach \cite{b013101,b013102,b013103,b013104,b102405,b013106,b013109,b111601} and the direct approach \cite{b110205,b110206,b013105,b013107,b013108,b020101}. The systematic study of the connections and differences between these two approaches has been conducted in \cite{b110205}.  The original problem of MFGs is that $N$-players are non-cooperative to optimize their own individual cost. The fixed-point approach takes the route from infinite-population to finite-population systems. This approach starts by constructing an auxiliary optimal control problem for a representative player through taking $N\rightarrow\infty$. Next, the fixed-point equation can be obtained by solving the auxiliary control problem together with the consistent condition.  Then, the decentralized control can be designed by solving the associated fixed-point equation. Finally, the designed decentralized control can be shown to be an $\varepsilon$-Nash equilibrium for the finite-population system. By this approach, decentralized Nash equilibria for various models can be effectively designed, but its challenge lies in verifying the solvability of the fixed-point equation. Therefore, merely the relatively conservative sufficient conditions for the solvability of fixed-point equations can usually be obtained \cite{b102405}. Conversely, the direct approach starts with the finite-population and extends to the infinite-population system. The first step is to solve the centralized Nash equilibrium of the original finite $N$-player Nash game under the global information pattern; the next step is to design the decentralized control by considering the limite of the centralized Nash equilibrium in the infinite-population system. This approach circumvents fixed-point equations, albeit at the expense of complex computations of centralized strategies \cite{b110206}. Not only the existence and uniqueness of centralized Nash equilibria are investigated by this approach, but the convergence issues are also examined \cite{b110203,b110202}. Therefore, this approach can obtain the necessary and sufficient conditions for the existence of decentralized control.   A detailed analysis can be found in \cite{b110205,b110206}.   Unlike the traditional decentralized control that only utilizes the individual' own state information, the direct approach is further used to design the asymptotically optimal distributed control according to a given network topology \cite{b111501}.

The open-loop and closed-loop strategies are important concepts in
control and game theory, and their relationship has been extensively studied by researchers in  single-player linear-quadratic (LQ) optimal control \cite{b102601} and two-player LQ differential games \cite{b102602,b102603}. It has been revealed that open-loop and closed-loop optimal controls are nearly equivalent for control problems, while open-loop and closed-loop Nash equilibria are significantly different for nonzero-sum differential games. The majority of MFG literature design decentralized control under the concept of closed-loop strategies \cite{b110202,b110204}, and some works are carried out  under the concept of open-loop strategies \cite{b110201,b110203}. Although the research routes and methods are various, an interesting observation is that under the different concepts of open-loop and closed-loop strategies, the designed asymptotically decentralized Nash equilibria for the same large-population system are identical \cite{b102201,b111701}. Naturally, we have the questions: Is this phenomenon universal? Why are the open-loop and closed-loop Nash equilibria significantly different for the finite-population system, but consistent for the infinite-population system? This inspires us to explore the impact of population scale on the design of open-loop and closed-loop strategies in LQ MFGs. 

In the first step of the fixed-point approach, an auxiliary optimal control problem is formulated within the infinite-population system. This auxiliary control problem is solved under the concept of open-loop \cite{b110301} and closed-loop \cite{b102405} strategies, respectively. Due to the equivalence of open-loop and closed-loop optimal control in control problems \cite{b102601}, the fixed-point approach cannot capture the impact of employing different open-loop and closed-loop strategies on the design of decentralized control. Conversely, in the first step of the direct approach, the centralized Nash equilibrium is solved for the finite-population system. Since there is a significant divergence between centralized open-loop and closed-loop Nash equilibria, the direct approach can better illustrate the impact of different strategies on the design of decentralized control. Therefore, in this paper we adopt the direct approach to investigate the impact of using different concept of open-loop and closed-loop strategies on the design of decentralized control strategies in MFGs. However, there is an absence of MFG literature employing a unified framework to concurrently investigate open-loop and closed-loop Nash equilibria in the finite-population system and designing decentralized control in the infinite-population system, resulting in challenges for comparative analysis. See \cite{b110205,b110206} for the closed-loop Nash equilibrium by an analytical method and \cite{b110201,b110204} for the open-loop Nash equilibrium by a probabilistic method, respectively.  Therefore, a unified framework is needed to  simultaneously study both open-loop and closed-loop strategies for  efficient analysis and comparison.  In \cite{b102603}, the open-loop and closed-loop Nash equilibria are studied for the two-person nonzero-sum differential games by using  forward-backward stochastic differential equations (FBSDEs) together with Riccati equations, and the differences and connections between these two strategies are analyzed in depth.  Inspired by the above works, we will study the centralized open-loop and closed-loop Nash equilibria for the finite-population system  using FBSDEs and Riccati equations within a unified framework in the first step of the direct approach. Thus, the impact of employing open-loop and closed-loop strategies on the design of decentralized control can be explored by examining the limit of centralized solutions in the infinite-population system. 

In this paper, we systematically study the effect of employing different concepts of open-loop and closed-loop strategies on the design of asymptotically decentralized Nash equilibria for LQ MFGs from finite to infinite-population systems. The large-population system contains $N$ homogeneous players, thus the system is symmetric and players are exchangeable. The state-weighting matrix in the cost functional is {indefinite}. In the first step of the direct approach, we solve the centralized open-loop and closed-loop Nash equilibria for the finite-population system under global information pattern. For the open-loop system, the necessary and sufficient conditions for the existence of centralized open-loop Nash equilibria are obtained by using the variational analysis, which are characterized by an adapted solution to a system of FBSDEs together with a convexity condition. Furthermore, we obtain the feedback representation of centralized open-loop Nash equilibrium by using the solution to a system of Riccati equations consisting of four ordinary differential equations (ODEs). For the closed-loop system, the necessary and sufficient conditions for the existence of centralized closed-loop Nash equilibria are obtained by the method of completing the square, which are characterized by the solvability of a system of Riccati equation  consisting of three ODEs. Since both centralized open-loop and closed-loop Nash equilibria are characterized by Riccati equations, they can be conveniently analyzed and compared. In the second step of the direct approach, we design the decentralized control under local information pattern by examining the limit of centralized results in the infinite-population system.  
The systems of Riccati equations in the finite-population system are studied for their limits as the population scale approaches infinity. The asymptotic solvability of the centralized open-loop and closed-loop Nash equilibria for the finite-population system is studied and is characterized by the solvability of a standard Riccati equation and an asymmetric Riccati equation. Subsequently, asymptotically decentralized Nash equilibria are designed by using mean field approximations. The main novelty of the paper are summarized as follows.
\begin{itemize}
	
	\item A unified framework is employed to study the open-loop and closed-loop Nash equilibria of LQ MFGs in both finite and infinite-population systems.   The impact of population scale on open-loop and closed-loop strategies is further elucidated through an analysis of the convergence issues of the associated systems of Riccati equations. It emerges that there is a clear distinction between centralized open-loop and closed-loop Nash equilibria in the finite-population system. However, when these strategies are applied to design asymptotically decentralized Nash equilibria in the infinite-population system, the outcomes coincide asymptotically.  This  indicates that the population scale plays a vital role in the design of open-loop and closed-loop strategies but the choice between open-loop and closed-loop strategies  is not crucial for  the design of decentralized control in LQ MFGs. 
	
	
	
	
	
	\

	\item The asymptotic solvability of the open-loop and closed-loop systems is investigated by examining the convergence issues of Riccati equations. The results show that the conditions for asymptotic solvability are identical for both cases, which are equivalent to the solvability conditions of a standard Riccati equation and an asymmetric Riccati equation. 
	
\end{itemize}

The rest of the paper is organized as follows. In Section \ref{formualation} we formulate the MFG problem, and introduce the  centralized open-loop and closed-loop Nash equilibria for the finite-population system, as well as the asymptotically decentralized Nash equilibria for the infinite-population system. 
In Section \ref{centralized}, the centralized open-loop and closed-loop Nash equilibria are obtained. The design of the asymptotically decentralized Nash equilibria is detailed in Section \ref{decentralized}. Finally, Section \ref{conclusion} concludes the paper.

Notation: Through the paper, $\mathbb{R}^{n\times m}$ denotes the set of all $n\times m$ real matrices. We use $|\cdot|$ and $\langle \cdot,\cdot\rangle$ to denote the Euclidean norm and Euclidean inner product. For a vector or matrix $M$, $M^T$ denotes its transpose. Consider a finite time horizon $[0,T]$. Suppose $(\Omega,\mathcal{F},\mathbb{F},\mathbb{P})$ is a complete filtered probability space. 
Let $L^p(0,T;\mathbb{R}^k)$ be the space of all $\mathbb{R}^k$-valued measurable functions $h(\cdot)$ on $[0,T]$ satisfying $\mathbb{E}\int_0^T|h(t)|^pdt<\infty$; $L^p_{\mathbb{F}}(0,T;\mathbb{R}^k)$ be the space of all $\mathbb{R}^k$-valued $\mathcal{F}_t$ progressively measurable processes $h(\cdot)$ on $[0,T]$ satisfying $\mathbb{E}\int_0^T|h(t)|^pdt<\infty$, where $1\leq p<\infty$.  For a set of vectors $\{x_1, x_2, . . .,x_N\}$, $\mathrm{vec}(x_1, x_2, . . .,x_N)$ denotes the vector $[x_1^T,x_2^T, . . .,x_N^T]^T$. For a set of matrices $\{A_1,A_2, . . .,A_N\}$, $\mathrm{cols}(A_1,A_2, . . .,A_N)$ denotes the matrices  $[A_1,A_2, . . .,A_N]$.

\section{Problem Formulation}\label{formualation}

\subsection{The $N$-player Nash game and centralized control}
Consider a large-population system with $N$ players denoted by $\mathcal{A}_i, i\in\mathcal{N}\triangleq\{1,2, . . .,N\}$. The dynamics of player $\mathcal{A}_i$, $i\in\mathcal{N}$, is given by the following stochastic differential equation
\begin{equation}\label{e100101}
	dx_i(t)=[Ax_i(t)+Gx^{(N)}(t)+Bu_i(t)]dt+\sigma dw_i(t),\quad x_i(0)=x_{i0},
\end{equation}
where $x_i(t)\in\mathbb{R}^n$ and $u_i(t)\in\mathbb{R}^m$ are the state and control of player $\mathcal{A}_i$, respectively. Players are coupled with each other through the mean field term $x^{(N)}(t)\triangleq (1/N)\sum_{i=1}^Nx_i(t)$. The coefficients of the state equation have compatible dimensions. $\{w_i(t),\forall i\in\mathcal{N}\}$ are a sequence of standard $1$-dimensional  Brownian motions  defined on a complete filtered probability space $(\Omega,\mathcal{F},\mathbb{F},\mathbb{P})$, where $\mathbb{F}\triangleq\{\mathcal{F}_t\}_{0\leq t\leq T}$ is the natural filtration of $\{w_i(t),\forall i\in\mathcal{N}\}$  augmented by all the $\mathbb{P}$-null sets in $\mathcal{F}$. For notational simplicity, we only consider constant coefficients for the model. Our analysis and results can be easily extended to the case of time-dependent coefficients.  Let $x^{(N)}_0\triangleq (1/N)\sum_{i=1}^Nx_{i0}$. We impose the following assumption:

\noindent\textbf{(A1)} The initial states  $\{x_{i0},$ $\forall i\in\mathcal{N}\}$ are mutually independent with $\mathbb{E}x_{i0}=\bar{x}_0$ and $\mathbb{E}|x_{i0}|^2\leq \infty$. The Brownian motions $\{w_i(t),\forall i\in\mathcal{N}\}$ are mutually independent and also independent of the initial states. 


For $i\in\mathcal{N}$, denote $\mathcal{N}_{-i}\triangleq \mathcal{N}-\{ i\}$, and $\mathbf{x}_0\triangleq \mathrm{vec}(x_{10},x_{20}, .. .,x_{N0})$. The cost functional of player $\mathcal{A}_i$ is given by 
\begin{equation}\label{e100102}
	\begin{split}		
		{\rm J}_i(\mathbf{x}_0;u_i(\cdot),u_j(\cdot),\forall j\in\mathcal{N}_{-i})= \mathbb{E}\int_0^T\Big[|x_i(t)-\Gamma x^{(N)}(t)-\eta|^2_{Q}+|u_i(t)|^2_{R}\Big]dt
		+ \mathbb{E}|x_i(T)-\Gamma_{f} x^{(N)}(T)-\eta_f|^2_{Q_{f}}. 
	\end{split}
\end{equation}%
The coefficients of the cost functional have compatible dimensions, where $Q$ and $Q_f$ are indefinite symmetric matrices, and $R$ is a positive-definite symmetric matrix. The admissible  centralized control set for player $\mathcal{A}_i,i\in\mathcal{N}$, is given by 
{ \begin{equation*}
		\begin{split}
			\mathcal{U}_c[0,T]\triangleq\Big\{u:[0,T]\times \Omega\rightarrow\mathbb{R}^{m}|\ &u(t)\in \sigma(x_{k0},w_k(s),s\leq t,\forall k\in\mathcal{N}),\mathbb{E}\int_0^T|u(s)|^2ds<\infty\Big\},
		\end{split} 
\end{equation*}}%
where the global information $\{x_{k0},w_k(t), \forall k\in\mathcal{N}\}$ of the large-population system is used.
Any element $u_i(\cdot)\in\mathcal{U}_c[0,T]$ is called an admissible centralized control of player $\mathcal{A}_i$ on $[0,T]$. In this setting, player $\mathcal{A}_i$ wishes to minimize the cost functional \eqref{e100102} by selecting a control process $u_i(\cdot)\in\mathcal{U}_c[0,T]$. We call the state equation \eqref{e100101} with the cost functional \eqref{e100102} the open-loop system of the large-population system in this paper.  For the centralized open-loop Nash equilibrium, we have the following definition:
\begin{definition}\label{def1}
	For the  open-loop system \eqref{e100101}-\eqref{e100102}, a  tuple $(\check{u}_i(\cdot)\in\mathcal{U}_c[0,T],$ $\forall i\in\mathcal{N})$ is called a centralized open-loop Nash equilibrium if for any $(u_i(\cdot)\in\mathcal{U}_c[0,T],$ $\forall i\in\mathcal{N})$, the following hold:
	\begin{equation}\label{e100205}
		\begin{split}
			{\rm J}_i(\mathbf{x}_0;\check{u}_i(\cdot),\check{u}_j(\cdot),\forall j\in\mathcal{N}_{-i})\leq {\rm J}_i(\mathbf{x}_0;u_i(\cdot),\check{u}_j(\cdot),\forall j\in\mathcal{N}_{-i}),\quad \forall i\in\mathcal{N}.
		\end{split}
	\end{equation}
\end{definition} 
If the centralized open-loop Nash equilibrium $(\check{u}_i(\cdot),\forall i\in\mathcal{N})$  exists and can be further represented as  
\begin{equation}\label{e041702}
	\begin{split}
		\check{u}_i(\cdot)=\sum_{l=1}^N\check{K}_{il}(\cdot)\check{x}_l(\cdot)+\check{v}_i (\cdot),\quad \forall i\in\mathcal{N},
	\end{split}
\end{equation}%
where $\check{K}_{il}(\cdot)\in  \mathcal{L}_1[0,T]\triangleq  L^2(0,T;\mathbb{R}^{m \times n})$ and  $\check{v}_i(\cdot)\in   \mathcal{L}_2[0,T]\triangleq  L^2(0,T;\mathbb{R}^{m \times 1})$, then we call 
$(\check{K}_{ij}(\cdot),\check{v}_i(\cdot),\forall i,j\in\mathcal{N})$ the feedback representation of the centralized open-loop Nash equilibrium $(\check{u}_i(\cdot),\forall i\in\mathcal{N})$. 

To further study the closed-loop control of this large-population system and compare it with open-loop control, we need to introduce the closed-loop system.   Consider the following closed-loop control 
\begin{equation} \label{e100701}
	u_i(t)=\sum_{l=1}^NK_{il}(t)x_l(t)+v_i(t),\quad i\in\mathcal{N},
\end{equation}
where $K_{il}(\cdot)\in   \mathcal{L}_1[0,T]$ and $v_i(\cdot)\in  \mathcal{L}_2[0,T]$.  Substituting the closed-loop control \eqref{e100701} into the state equation \eqref{e100101}, we obtain the following dynamics for player $\mathcal{A}_i$:
\begin{equation}\label{e100702}
	\begin{split}
		dx_i(t)=[Ax_i(t)\!+\!Gx^{(N)}(t)\!+\!\sum_{l=1}^NBK_{il}(t)x_l(t)\!+\!Bv_i(t)]dt\!+\! \sigma dw_i(t),\quad \forall i\in\mathcal{N}.
	\end{split}
\end{equation}%
With the solution $\mathbf{x}(\cdot)\triangleq \mathrm{vec}(x_1(\cdot),x_2(\cdot), . . . ,x_N(\cdot))$ to \eqref{e100702}, substituting the closed-loop control \eqref{e100701} into the cost functional \eqref{e100102}, we obtain 
\begin{equation}\label{e100703}
	\begin{split}		
		&{\rm J}_i(\mathbf{x}_0;K_i(\cdot)\mathbf{x}(\cdot)+v_i(\cdot),K_j(\cdot)\mathbf{x}(\cdot)+v_j(\cdot),\forall j\in\mathcal{N}_{-i})\\=\ &\mathbb{E}\int_0^T \Big[|x_i(t)-\Gamma x^{(N)}(t)-\eta|^2_{Q}+|	\sum_{l=1}^NK_{il}(t)x_l(t)+v_i(t)|^2_{R}\Big]dt+\mathbb{E}|x_i(T)-\Gamma_{f} x^{(N)}(T)-\eta_f|^2_{Q_{f}},
	\end{split}
\end{equation}
where $K_i(\cdot)\triangleq \mathrm{cols}(K_{i1}(\cdot),K_{i2}(\cdot), . . .,K_{iN}(\cdot))$.  Similarly, we  can define ${\rm J}_i(\mathbf{x}_0;u_i(\cdot),$ $K_j(\cdot)\mathbf{x}(\cdot)+v_j(\cdot),\forall j\in\mathcal{N}_{-i})$ and ${\rm J}_i(\mathbf{x}_0;K_i(\cdot)\mathbf{x}(\cdot)+v_i(\cdot),u_j(\cdot),\forall j\in\mathcal{N}_{-i})$.  We call  $(u_i(\cdot),\forall i\in\mathcal{N})$  with $u_i(\cdot)$ defined by \eqref{e100701} the closed-loop control, and also call the tuple $(K_{ij}(\cdot),v_i(\cdot),\forall i,j\in\mathcal{N})$ the closed-loop control as long as no confusion arises. 
We call the state equation \eqref{e100702} with the cost functional \eqref{e100703} the closed-loop system of the large-population system under the closed-loop control $(K_{ij}(\cdot),v_i(\cdot),\forall i,j\in\mathcal{N})$. 
We now introduce the following definition of the centralized closed-loop Nash equilibrium.
\begin{definition} 	For the closed-loop system \eqref{e100702}-\eqref{e100703},  a tuple $(\hat{K}_{ij}(\cdot)\in \mathcal{L}_1,$ $\hat{v}_i(\cdot)\in \mathcal{L}_2[0,T],\forall i,j\in\mathcal{N})$
	is called a centralized closed-loop Nash equilibrium  if for any  $(K_{ij}(\cdot)\in \mathcal{L}_1[0,T],v_i(\cdot),\in \mathcal{L}_2[0,T],\forall i,j\in\mathcal{N})$ the following hold:
	\begin{equation}\label{e100206b}
		\begin{split}
			{\rm J}_i(\mathbf{x}_0;\hat{K}_i(\cdot)\hat{\mathbf{x}}(\cdot)+\hat{v}_i(\cdot),\hat{K}_j(\cdot)\hat{\mathbf{x}}(\cdot)+\hat{v}_j(\cdot),\forall j\in\mathcal{N}_{-i})\leq {\rm J}_i(\mathbf{x}_0;u_i(\cdot),\hat{K}_j(\cdot)\mathbf{x}(\cdot)+\hat{v}_j(\cdot),\forall j\in\mathcal{N}_{-i}),\quad \forall i\in\mathcal{N}.
		\end{split}
	\end{equation}
\end{definition}

As the first step of the direct approach, we will solve the centralized open-loop and closed-loop Nash equilibria of the $N$-player system and analyze their differences.  The main problems are described as follows.

\noindent\textbf{Problem (P1).} For the $N$-player open-loop system \eqref{e100101}-\eqref{e100102}, seek the feedback representation $(\check{K}_{ij}(\cdot),$ $\check{v}_i(\cdot),\forall i,j\in\mathcal{N})$ of the centralized open-loop Nash equilibrium $(\check{u}_i(\cdot),\forall i\in\mathcal{N})$  satisfying \eqref{e100205}.

\noindent\textbf{Problem (P2).} For the $N$-player closed-loop system \eqref{e100702}-\eqref{e100703}, seek the centralized closed-loop Nash equilibrium $(\hat{K}_{ij}(\cdot),\hat{v}_i(\cdot),\forall i,j\in\mathcal{N})$ satisfying \eqref{e100206b}.

\subsection{Decentralized control}




It is not realistic for players to collect the global information of the large-population system to implement centralized control strategies due to the limitations of sensing, computing and communication capabilities. Alternatively, a more practical approach is for each player to utilize their own state information to design decentralized control.  The admissible decentralized control set for $\mathcal{A}_i$, $i\in\mathcal{N}$, is given as follows.
\begin{equation*}
	\begin{split}
		\mathcal{U}_d^i[0,T] \triangleq  \Big\{u:[0,T]\times \Omega\rightarrow\mathbb{R}^{m}\ |\ u(t)\in \sigma(x_i(s),s\leq t),\mathbb{E}\int_0^T|u(s)|^2ds<\infty\Big\}. 
	\end{split}
\end{equation*}%
Any element $u_i(\cdot)\in\mathcal{U}_d^i[0,T]$ is called an admissible decentralized control of player $\mathcal{A}_i$ on $[0,T]$. Attaining a Nash equilibrium is challenging with merely local information. Consequently, our goal is to design an approximate Nash equilibrium as closely as possible, namely the decentralized $\varepsilon$-Nash equilibrium, which is defined as follows.
\begin{definition}
	For a given $\varepsilon\geq 0$, a $N$-tuple $(u_i^d(\cdot)\in\mathcal{U}^i_d[0,T],\forall i\in\mathcal{N})$ is called a decentralized $\varepsilon$-Nash equilibrium if the following hold for any $i\in\mathcal{N}$:
	\begin{equation}\label{e101001}
		\begin{split}
			&{\rm J}_i(\mathbf{x}_0;u_i^d(\cdot),u_j^d(\cdot),j\in\mathcal{N}_{-i}) \leq  \inf_{u_i(\cdot)\in\mathcal{U}_c[0,T]}{\rm J}_i(\mathbf{x}_0;u_i(\cdot),u_j^d(\cdot),j\in\mathcal{N}_{-i})+\varepsilon.
		\end{split}
	\end{equation}
\end{definition}
For the decentralized control, our aim is not only to design a decentralized $\varepsilon$-Nash equilibrium for a fixed $\varepsilon > 0$ but also to have $\varepsilon \rightarrow 0$ as $N \rightarrow \infty$. That is, the designed decentralized control can achieve almost  the same effect as the centralized Nash equilibrium in the large-population system. Therefore we define the asymptotically decentralized Nash equilibrium as follows. 
\begin{definition}
	A $N$-tuple $(u_i^d(\cdot)\in\mathcal{U}_d^i[0,T],\forall i\in\mathcal{N})$ is called an asymptotically decentralized Nash equilibrium with respect to the number  of players in the large-population if the following hold for any $i\in\mathcal{N}$:
	{ \setlength{\abovedisplayskip}{4pt}
		\setlength{\belowdisplayskip}{4pt}
		\setlength{\abovedisplayshortskip}{4pt}
		\setlength{\belowdisplayshortskip}{4pt}\begin{equation}\label{e100305}
			\begin{split}
				{\rm J}_i(\mathbf{x}_0;u_i^d(\cdot),u_j^d(\cdot),j\in\mathcal{N}_{-i}) \leq \! \inf_{u_i(\cdot)\in\mathcal{U}_c[0,T]} {\rm J}_i(\mathbf{x}_0;u_i(\cdot),u_j^d(\cdot),j\in\mathcal{N}_{-i})+O(\frac{1}{N}).\!\!
			\end{split}
	\end{equation}}
\end{definition}

As the second step of the direct approach, we will design the asymptotically  decentralized Nash equilibria. The main problems are stated below. 

\noindent\textbf{Problem (P3).} Based on the feedback representation $(\check{K}_{ij}(\cdot),\check{v}_i(\cdot),\forall i,j\in\mathcal{N})$ of the centralized open-loop Nash equilibrium $(\check{u}_i(\cdot),\forall i\in\mathcal{N})$, design an asymptotically decentralized Nash equilibrium $(u_i^d(\cdot),\forall i\in\mathcal{N})$  
satisfying \eqref{e100305}.

\noindent\textbf{Problem (P4).} Based on the centralized closed-loop Nash equilibrium $(\hat{K}_{ij}(\cdot),\hat{v}_i(\cdot),$ $\forall i,j\in\mathcal{N})$, design an asymptotically decentralized Nash equilibrium $(u_i^d(\cdot),\forall i\in\mathcal{N})$  
satisfying \eqref{e100305}.

\section{Main Results for Centralized Control}\label{centralized}
In this section we study the centralized open-loop and closed-loop Nash equilibria for the finite-population system within a unified technical framework. Centralized open-loop and closed-loop Nash equilibria characterized by FBSDEs can be obtained through variational analysis. By decoupling the FBSDEs with Riccati equations, explicit solutions can  be further derived, which are beneficial for subsequent analysis and comparison.




\subsection{Centralized open-loop Nash equilibrium to Problem (P1)}\label{subsection3.1}
By using the variational analysis, the main results for the existence of centralized open-loop Nash equilibria can be obtained as follows.
\begin{theorem}\label{th1}
	Let   $(\check{x}_i(\cdot),\forall i\in\mathcal{N})$ and $\check{x}^{(N)}(\cdot)$ are respectively the corresponding states and the mean field term under the control $(\check{u}_i(\cdot)\in\mathcal{U}_c[0,T],\forall i\in\mathcal{N})$. Then $(\check{u}_i(\cdot),\forall i\in\mathcal{N})$ is a centralized open-loop Nash equilibrium of Problem (P1) if and only if the following system of FBSDEs 
	{ \begin{equation}\label{e101101}
			\begin{split}
				\left\{
				\begin{aligned}
					d\check{x}_i(t)=&\ [A\check{x}_i(t)+G\check{x}^{(N)}(t)-BR^{-1} B^T\check{p}_i^i(t)]dt+\sigma dw_i(t),\quad \check{x}_i(0)=x_{i0}, \\
					d\check{p}_i^i(t)=&-\Big[A^T\check{p}_i^i(t)+G^T\check{p}^{i,(N)}(t)+(Q-\frac{\Gamma^TQ}{N})\check{x}_i(t)
					-(Q\Gamma-\frac{\Gamma^TQ\Gamma}{N}) \check{x}^{(N)}(t)-(Q-\frac{\Gamma^TQ}{N})\eta\Big]dt\\&+\sum_{k=1}^N\check{z}_{ik}^i(t)dw_k(t),\quad 
					\check{p}_i^i(T)= (Q_f\!-\!\frac{\Gamma^T_fQ_f}{N})\check{x}_i(T)-(Q_f\Gamma_f\!-\!\frac{\Gamma^T_fQ_f\Gamma_f}{N}) \check{x}^{(N)}(T)-(Q_f-\frac{\Gamma^T_fQ_f}{N})\eta_f,\\ 
					d\check{p}_j^i(t)=&-\Big[A^T\check{p}_j^i(t)+G^T\check{p}^{i,(N)}(t)-\frac{\Gamma^TQ}{N}\check{x}_i(t)+\frac{\Gamma^TQ\Gamma}{N}\check{x}^{(N)}(t)+\frac{\Gamma^TQ}{N}\eta\Big]dt\\&+\sum_{k=1}^N\check{z}_{jk}^i(t)dw_k(t),\quad
					\check{p}_j^i(T)=-\frac{\Gamma^T_fQ_f}{N}\check{x}_i(T)+\frac{\Gamma^T_fQ_f\Gamma_f}{N}\check{x}^{(N)}(T)+\frac{\Gamma^T_fQ_f}{N}\eta_f,\ \forall i\in\mathcal{N},\ \forall j\in\mathcal{N}_{-i}, 				
				\end{aligned}
				\right.
			\end{split}
	\end{equation}}%
				admits an $\mathcal{F}_t$-adapted solution $(\check{x}_i(t),\check{p}_j^i(t),\check{z}^i_{jk}(t),\forall i,j,k\in\mathcal{N})$ on $t\in[0,T]$, where $\check{p}^{i,(N)}(t)\triangleq (1/N)\sum_{j=1}^N\check{p}^i_j(t)$, and the following convexity  condition holds $\forall u_i(\cdot)\in\mathcal{U}_c[0,T]$ and $\forall i\in\mathcal{N}$
				\begin{equation}\label{e091402}
					\mathbb{E}\int_0^T\Big[|\check{y}_i^i(t)-\Gamma \check{y}^{i,(N)}(t)|^2_Q+|u_i(t)|^2_R\Big]dt+\mathbb{E}|\check{y}_i^i(T)-\Gamma_f \check{y}^{i,(N)}(T)|^2_{Q_f}\geq 0
				\end{equation}
				where $\check{y}^{i,(N)}(t)=(1/N)\sum_{k=1}^N\check{y}^i_k(t)$ and $\check{y}^i_k(t), \forall k\in\mathcal{N}$, satisfy the following ODEs
				\begin{subequations}\label{e101305}
					\begin{align}
						\dot{\check{y}}_i^i(t)=&\ A\check{y}_i^i(t)+G\check{y}^{i,(N)}(t)+Bu_i(t),\quad \check{y}_i^i(0)=0,\\
						\dot{\check{y}}_j^i(t)=&\ A\check{y}_j^i(t)+G\check{y}^{i,(N)}(t),\quad \check{y}_j^i(0)=0,\quad \forall j\in\mathcal{N}_{-i}.
					\end{align}
				\end{subequations}	
				In this case, the centralized open-loop Nash equilibrium  is given by 
				\begin{equation}\label{e101102}
					\check{u}_i(t)=-R^{-1} B^T\check{p}_i^i(t),\quad \forall i\in\mathcal{N}.
				\end{equation}
				
			\end{theorem}
			
			\textbf{Proof}   See Appendix \ref{proofofthe1}. $\hfill\blacksquare$
			
			Inspired by the decoupling technique used in \cite{b020401,b110703,b110601}, we now look at the solvability of FBSDEs \eqref{e101101}.  Denote $\Upsilon\triangleq BR^{-1}B^T$.  We introduce the following system of Riccati equations:
			\begin{subequations}\label{e101201}
				\begin{align}
					\begin{split}\label{e101201a}
						\left\{
						\begin{aligned}
							&\dot{\check{P}}^N_1(t)=\check{P}^N_1(t)\Upsilon \check{P}^N_1(t)+(N-1)\check{P}^N_2(t)\Upsilon \check{P}^N_2(t)-\check{P}^N_1(t)(A+\frac{G}{N})-(A+\frac{G}{N})^T\check{P}^N_1(t)\\&\quad\quad\quad\quad-\check{P}^N_2(t)\frac{N-1}{N}G-\frac{N-1}{N}G^T\check{P}^N_4(t)-Q+\frac{\Gamma^TQ}{N}+\frac{Q\Gamma}{N}-\frac{\Gamma^TQ\Gamma}{N^2},\\& 	\check{P}^N_1(T)=Q_f-\frac{\Gamma^T_fQ_f}{N}-\frac{Q_f\Gamma_f}{N}+\frac{\Gamma^T_fQ_f\Gamma_f}{N^2},
						\end{aligned}
						\right.
					\end{split}\\
					\begin{split}\label{e101201b}
						\left\{
						\begin{aligned}
							&\dot{\check{P}}^N_2(t)=\check{P}^N_1(t)\Upsilon \check{P}^N_2(t)+\check{P}^N_2(t)\Upsilon \check{P}^N_1(t)+(N-2)\check{P}^N_2(t)\Upsilon \check{P}^N_2(t)-\check{P}^N_1(t)\frac{G}{N}-(A+\frac{G}{N})^T\check{P}^N_2(t)\\&\quad\quad\quad\quad-\check{P}^N_2(t)(A+\frac{N-1}{N}G)-\frac{N-1}{N}G^T\check{P}^N_3(t)+\frac{Q\Gamma}{N}-\frac{\Gamma^TQ\Gamma}{N^2},\\& 	\check{P}^N_2(T)=\frac{\Gamma^T_fQ_f\Gamma_f}{N^2}-\frac{Q_f\Gamma_f}{N},
						\end{aligned}
						\right.
					\end{split}\\
					\begin{split}\label{e101201c}
						\left\{
						\begin{aligned}
							&\dot{\check{P}}^N_3(t)=\check{P}^N_4(t)\Upsilon \check{P}^N_2(t)+\check{P}^N_3(t)\Upsilon \check{P}^N_1(t)+(N-2)\check{P}^N_3(t)\Upsilon \check{P}^N_2(t)-\check{P}^N_3(t)(A+\frac{N-1}{N}G)\\&\quad\quad\quad\quad-\check{P}^N_4(t)\frac{G}{N}-(A+\frac{N-1}{N}G)^T\check{P}^N_3(t)-\frac{G^T}{N}\check{P}^N_2(t)-\frac{\Gamma^TQ\Gamma}{N^2},\quad 	\\&\check{P}^N_3(T)=\frac{\Gamma^T_fQ_f\Gamma_f}{N^2},
						\end{aligned}
						\right.
					\end{split}\\
					\begin{split}\label{e101201d}
						\left\{
						\begin{aligned}
							&\dot{\check{P}}^N_4(t)=\check{P}^N_4(t)\Upsilon \check{P}^N_1(t)+(N\!-\!1)\check{P}^N_3(t)\Upsilon \check{P}^N_2(t)-\check{P}^N_4(t)(A+\frac{G}{N})+\frac{\Gamma^TQ}{N}\\&\quad\quad\quad\quad-	(A+\frac{N\!-\!1}{N}G)^T	\check{P}^N_4(t)-\frac{N\!-\!1}{N}\check{P}^N_3(t)G-\frac{G^T}{N}\check{P}^N_1(t)-\frac{\Gamma^TQ\Gamma}{N^2},	\\&\check{P}^N_4(T)=-\frac{\Gamma^T_fQ_f}{N}+\frac{\Gamma^T_fQ_f\Gamma_f}{N^2},
						\end{aligned}
						\right.
					\end{split}
				\end{align}	
			\end{subequations}
			and the following system of ODEs:
			\begin{subequations} \label{e042401}
				\begin{align}
					\begin{split}
						\left\{
						\begin{aligned}
							&\dot{\check{s}}^N_1(t)=-\Big[A+\frac{G}{N}-\Upsilon\check{P}^N_1(t)^T-(N-1)\Upsilon\check{P}^N_2(t)^T\Big]^T\check{s}^N_1(t)-\frac{N-1}{N}G^T\check{s}^N_2(t)+(Q-\frac{\Gamma^TQ}{N})\eta,\quad \\&\check{s}^N_1(T)=(\frac{\Gamma_f^TQ_f}{N}-Q_f)\eta_f,
						\end{aligned}
						\right.
					\end{split}\\
					\begin{split}
						\left\{
						\begin{aligned}
							&\dot{\check{s}}^N_2(t)=-\Big[\frac{G}{N}-\Upsilon\check{P}^N_4(t)^T-(N-1)\Upsilon\check{P}^N_3(t)^T\Big]^T\check{s}^N_1(t) -(A+\frac{N-1}{N}G)^T\check{s}^N_2(t)-\frac{\Gamma^TQ}{N}\eta,\quad\\ 	&\check{s}^N_2(T)=\frac{\Gamma_f^TQ_f}{N}\eta_f. 
						\end{aligned} 
						\right.
					\end{split}
				\end{align}	 
			\end{subequations}
			Then we can further obtain the following results about the feedback representation of centralized open-loop Nash equilibria.
			\begin{theorem}\label{closedlooppro}
				If the convexity condition \eqref{e091402} holds and the system of Riccati equations \eqref{e101201} admit a solution $(\check{P}^N_i(\cdot),\forall i=1,. . .,4)$, then the centralized open-loop Nash equilibrium $(\check{u}_i(\cdot),\forall i\in\mathcal{N})$ of Problem (P1) has the feedback representation 
				{ \setlength{\abovedisplayskip}{4pt}
					\setlength{\belowdisplayskip}{-1pt}
					\setlength{\abovedisplayshortskip}{4pt}
					\setlength{\belowdisplayshortskip}{-1pt}
					\begin{equation} \label{e100804}
						\begin{split}
							\check{u}_i(t)=&\ \check{K}^N_1(t)\check{x}_i(t)+\sum_{j\in\mathcal{N}_{-i}}\check{K}^N_2(t)\check{x}_j(t)+\check{v}^N(t),\quad \forall i\in\mathcal{N},
						\end{split} 
				\end{equation}}
				with 		
				\begin{align}\label{e122403a}
					\check{K}^N_1(t)= -R^{-1} B^T\check{P}_1^N(t),\ 
					\check{K}^N_2(t)=-R^{-1} B^T\check{P}_2^N(t),\  \check{v}^N(t)=	-R^{-1} B^T\check{s}_1^N(t). 
				\end{align}
				
				%
				%
			\end{theorem}

			\textbf{Proof}  
			See Appendix \ref{proofofclosedlooppro}.$\hfill\blacksquare$

			\subsection{Centralized closed-loop Nash equilibrium to Problem (P2)}\label{subsection3.2} In this subsection, we study Problem (P2). We first give the equivalent definition of the closed-loop Nash equilibria, thereby establishing the connection between open-loop and closed-loop systems, which plays an important role in studying the necessity conditions for the existence of closed-loop Nash equilibria. Some similar results of the connection between open-loop and closed-loop strategies can be found in \cite[Proposition 3.3]{b102601} and \cite[Proposition 3.2]{b110701},  so we omit the proof here. 
			
			\begin{lemma}\label{prop1}
				Let $\hat{K}_{ij}(\cdot)\in \mathcal{L}_1[0,T], \hat{v}_i(\cdot)\in \mathcal{L}_2[0,T]$ and $\hat{K}_i(\cdot)\triangleq \mathrm{cols}(\hat{K}_{i1}(\cdot),\hat{K}_{i2}(\cdot),$ $ . . .,\hat{K}_{iN}(\cdot)),\forall i,j\in\mathcal{N}$, then the following statements are equivalent:
				\begin{enumerate}
					\item The tuple $(\hat{K}_i(\cdot),\hat{v}_i(\cdot),\forall i\in\mathcal{N})$
					is called a centralized closed-loop Nash equilibrium for the closed-loop system \eqref{e100702}-\eqref{e100703}.\label{223}
					\item \label{9291} For any  $K_{ij}(\cdot)\in \mathcal{L}_1[0,T],v_i(\cdot)\in \mathcal{L}_2[0,T],\forall i,j\in\mathcal{N},$ the following hold:
					\begin{equation}\label{e101401}
						\begin{split}
							{\rm J}_i(\mathbf{x}_0;\hat{K}_i(\cdot)\hat{\mathbf{x}}(\cdot)\!+\!\hat{v}_i(\cdot),\hat{K}_j(\cdot)\hat{\mathbf{x}}(\cdot)\!+\!\hat{v}_j(\cdot),\forall j\in\mathcal{N}_{-i})\leq {\rm J}_i(\mathbf{x}_0;K_i(\cdot)\mathbf{x}(\cdot)\!+\!v_i(\cdot),\hat{K}_j(\cdot)\mathbf{x}(\cdot)\!+\!\hat{v}_j(\cdot),\forall j\in\mathcal{N}_{-i}).
						\end{split}
					\end{equation}
					
					\item For any  $v_i(\cdot)\in \mathcal{L}_2[0,T]$ 
					the following hold:
					\begin{equation}\label{e101402}
						\begin{split}
							{\rm J}_i(\mathbf{x}_0;\hat{K}_i(\cdot)\hat{\mathbf{x}}(\cdot)\!+\!\hat{v}_i(\cdot),\hat{K}_j(\cdot)\hat{\mathbf{x}}(\cdot)\!+\!\hat{v}_j(\cdot),\forall j\in\mathcal{N}_{-i})\leq {\rm J}_i(\mathbf{x}_0;\hat{K}_i(\cdot)\mathbf{x}(\cdot)\!+\!v_i(\cdot),\hat{K}_j(\cdot)\mathbf{x}(\cdot)\!+\!\hat{v}_j(\cdot),\forall j\in\mathcal{N}_{-i}).
						\end{split}
					\end{equation}
				\end{enumerate}
			\end{lemma}
			
			
			
			Next we consider the state equation 
			{\setlength{\abovedisplayskip}{4pt}
				\setlength{\belowdisplayskip}{4pt}
				\setlength{\abovedisplayshortskip}{4pt}
				\setlength{\belowdisplayshortskip}{4pt}
				\begin{align}\label{e101403}
					\begin{split}
						dx_i(t)=  [Ax_i(t)+Gx^{(N)}(t)+\sum_{l=1}^NB\hat{K}_{il}(t)x_l(t)+Bv_i(t)]dt + \sigma dw_i(t), 
					\end{split}
			\end{align}}
			with $x_i(0)=x_{i0}$, $\forall i\in\mathcal{N}$, and the cost functional  
			{\setlength{\abovedisplayskip}{4pt}
				\setlength{\belowdisplayskip}{4pt}
				\setlength{\abovedisplayshortskip}{4pt}
				\setlength{\belowdisplayshortskip}{4pt}
				\begin{equation}\label{e101404}
					\begin{split}		
						&\ {\rm \hat{J}}_i(\mathbf{x}_0;v_i(\cdot),v_{-i}(\cdot))\\\triangleq &\ \mathbb{E}\int_0^T \Big[|x_i(t)-\Gamma x^{(N)}(t)-\eta|^2_{Q}+|\sum_{l=1}^N\hat{K}_{il}(t)x_l(t)+v_i(t)|^2_{R}\Big]dt+\mathbb{E}|x_i(T)-\Gamma_{f} x^{(N)}(T)-\eta_f|^2_{Q_{f}}. 
					\end{split}
			\end{equation}}
			Then  we can construct the following auxiliary $N$-player open-loop Nash game: 
			
			\noindent\textbf{Problem (P5).} For the open-loop system \eqref{e101403}-\eqref{e101404}, seek a centralized open-loop Nash equilibrium $(\hat{v}_i(\cdot)\in   \mathcal{U}_c[0,T],$ $\forall i\in\mathcal{N})$
			for the initial $\mathbf{x}_{0}$.
			
			It is shown by 3) of Lemma \ref{prop1} that if $(\hat{K}_{ij}(\cdot),\hat{v}_i(\cdot),\forall i,j\in\mathcal{N})$ is a closed-loop Nash equilibrium for the closed-loop system \eqref{e100702}-\eqref{e100703}, then Problem (P5) is solvable and $(\hat{v}_i(\cdot),\forall i\in\mathcal{N})$ is the centralized open-loop Nash equilibrium for the open-loop system \eqref{e101403}-\eqref{e101404}. Therefore we can use the same procedure in Theorem \ref{th1} to solve Problem (P5). The results are stated as follows.  
			\begin{proposition}\label{prop2}
				Let $(\hat{K}_{ij}(\cdot)\in \mathcal{L}_1[0,T],\hat{v}_i(\cdot)\in \mathcal{L}_2[0,T],\forall i,j\in\mathcal{N})$ be a centralized closed-loop Nash equilibrium of Problem (P2). Then for any $\mathbf{x}_{0}\in \mathbb{R}^{Nn}$, the following system of FBSDEs 
				{ 
					\begin{equation}\label{e23101001}
						\begin{split}
							\left\{
							\begin{aligned}
								d\hat{x}_i(t)=&\ [A\hat{x}_i(t)+G\hat{x}^{(N)}(t)+\sum_{l=1}^NB\hat{K}_{il}(t)\hat{x}_l(t)+B\hat{v}_i(t)]dt+ \sigma dw_i(t), \quad \hat{x}_i(0)= x_{i0},\\
								d\hat{p}_i^i(t)=&  -[A^T\hat{p}^i_i(t)+G^T\hat{p}^{i,(N)}(t)+\sum_{l=1}^N\hat{K}_{li}^T(t)B^T\hat{p}^i_l(t) + 	\sum_{l=1}^N\hat{K}^T_{ii}(t)R\hat{K}_{il}(t)\hat{x}_l(t)+\hat{K}_{ii}^T(t)R \hat{v}_i(t) \\&+(Q-\frac{\Gamma^TQ}{N})\hat{x}_i(t)  -(Q\Gamma-\frac{\Gamma^TQ\Gamma}{N}) \hat{x}^{(N)}(t)-(Q-\frac{\Gamma^TQ}{N})\eta]dt+\sum_{l=1}^N\hat{z}^i_{il}(t)dw_l(t),  \\
								\hat{p}^i_i(T)=&\ (Q_f\!-\!\frac{\Gamma_f^TQ_f}{N})\hat{x}_i(T)\!-\!(Q_f\Gamma_f\!+\!\frac{\Gamma^T_fQ_f\Gamma_f}{N}) \hat{x}^{(N)}(T)\!+\!(\frac{\Gamma^T_fQ_f}{N} \!-\!Q_f)\eta_f,\\
								d\hat{p}_i^j(t)=&  -[A^T\hat{p}^i_j(t)+G^T\hat{p}^{i,(N)}(t)+\sum_{l=1}^N\hat{K}_{lj}^T(t)B^T\hat{p}^i_l(t)+\hat{K}_{ij}^T(t)R \hat{v}_i(t) + 	\sum_{l=1}^N\hat{K}^T_{ij}(t)R\hat{K}_{il}(t)\hat{x}_l(t)\\& -\frac{\Gamma^TQ}{N}\hat{x}_i(t) +\frac{\Gamma^TQ\Gamma}{N} \hat{x}^{(N)}(t) +\frac{\Gamma^TQ}{N}\eta]dt+\sum_{l=1}^N\hat{z}^i_{jl}(t)dw_l(t),\\
								\hat{p}^i_j(T)=&\ -\frac{\Gamma^T_fQ_f}{N}\hat{x}_i(T)+\frac{\Gamma^T_fQ_f\Gamma_f}{N} \hat{x}^{(N)}(T)+\frac{\Gamma_f^TQ_f}{N}\eta_f,\  \forall i\in\mathcal{N},\  \forall j\in\mathcal{N}_{-i},
							\end{aligned}
							\right.
						\end{split}
				\end{equation}}%
				admits an $\mathcal{F}_t$-adapted solution $(\hat{x}_i,\hat{p}^i_j,\hat{z}^i_{jl},\forall i,j,l\in\mathcal{N})$ and the following stationarity condition holds: 
				\begin{equation}\label{e23101603}
					B^T\hat{p}^i_i(t)+	\sum_{l=1}^NR\hat{K}_{il}(t)\hat{x}_l(t)+R\hat{v}_i(t)=0.
				\end{equation}
			\end{proposition}
			\textbf{Proof}  
			According to  Lemma \ref{prop1},   $(\hat{K}_{ij}(\cdot),\hat{v}_i(\cdot),\forall i,j\in\mathcal{N})$ is a centralized closed-loop Nash equilibrium for the closed-loop system \eqref{e100702}-\eqref{e100703} is equivalent to that $(\hat{v}_i(\cdot),\forall i\in\mathcal{N})$ is an open-loop Nash equilibrium for the open-loop system \eqref{e101403}-\eqref{e101404}. We can solve Problem (P5) by the variational method used in Theorem \ref{th1} and obtain the necessary condition for the existence of $(\hat{v}_i(\cdot),\forall i\in\mathcal{N})$. Comparing the cost functional \eqref{e101404} with \eqref{e100102}, the open-loop system \eqref{e101403}-\eqref{e101404} with respect to the control $(v_i(\cdot),\forall i\in\mathcal{N})$ is more complicated than the open-loop system \eqref{e100101}-\eqref{e100102} with respect to the control $(u_i(\cdot),\forall i\in\mathcal{N})$ due to the intersection terms of state and control in \eqref{e101404}. This proposition cannot simply obtained by using variable substitutions in Theorem \ref{th1}. But the steps of using the variational method are the same, the system of FBSDEs \eqref{e23101001} for the existence of open-loop Nash equilibrium and the stationarity condition \eqref{e23101603} satisfied by $(\hat{K}_{ij}(\cdot),\hat{v}_i(\cdot),\forall i,j\in\mathcal{N})$ can be similarly obtained.$\hfill\blacksquare$

					Because of the homogeneity of players, the large-population system is symmetric and players are exchangeable. Therefore the closed-loop controllers  should have the same structure. Also inspired by the feedback representation of the  centralized open-loop Nash equilibria \eqref{e100804}, we consider the closed-loop control \eqref{e100701} with the following structure
					\begin{equation} \label{e23101701}
						K_{ii}(t)=K^N_1(t),\quad K_{ij}(t)=K^N_2(t),\quad v_i(t)=v^N(t),\quad \forall i\in\mathcal{N},\quad \forall j\in\mathcal{N}_{-i}.
					\end{equation} 
					

					By using the symmetric property of the large-population system, the system of FBSDEs \eqref{e23101001} can be rewritten as  
					{	\begin{equation}\label{e23101703}
							\begin{split}
								\left\{
								\begin{aligned}
									d\hat{x}_i(t)=&\ \Big\{[A+B\hat{K}^N_1(t)-B\hat{K}^N_2(t)]\hat{x}_i(t)+[G+NB\hat{K}^N_2(t)]\hat{x}^{(N)}(t) +B\hat{v}^N(t)\Big\}dt+ \sigma dw_i(t), \quad \hat{x}_i(0)= x_{i0}\\
									d\hat{p}_i^i(t)=&  -\Big\{[A+B\hat{K}^N_1(t)-B\hat{K}^N_2(t)]^T\hat{p}^i_i(t)+[G+NB\hat{K}^N_2(t)]^T\hat{p}^{i,(N)}(t) + [\hat{K}^N_1(t)^TR\hat{K}^N_1(t)\\&-\hat{K}^N_1(t)^TR\hat{K}^N_2(t)+Q-\frac{\Gamma^TQ}{N}]	\hat{x}_i(t)+[N\hat{K}^N_1(t)^TR\hat{K}^N_2(t)-Q\Gamma+\frac{\Gamma^TQ\Gamma}{N}]\hat{x}^{(N)}(t)\\&+\hat{K}^N_1(t)^TR \hat{v}^N(t)  -(Q-\frac{\Gamma^TQ}{N})\eta\Big\}dt+\sum_{l=1}^N\hat{z}^i_{il}(t)dw_l(t), \\
									\hat{p}^i_i(T)=&\ (Q_f\!-\!\frac{\Gamma_f^TQ_f}{N})\hat{x}_i(T)\!-\!(Q_f\Gamma_f\!+\!\frac{\Gamma^T_fQ_f\Gamma_f}{N}) \hat{x}^{(N)}(T)\!+\!(\frac{\Gamma^T_fQ_f}{N} \!-\!Q_f)\eta_f,\\
									d\hat{p}^i_j(t)=&  -\Big\{[A+B\hat{K}^N_1(t)-B\hat{K}^N_2(t)]^T\hat{p}^i_j(t)+[G+NB\hat{K}^N_2(t)]^T\hat{p}^{i,(N)}(t)+ [\hat{K}^N_2(t)^TR\hat{K}^N_1(t)\\&-\hat{K}^N_2(t)^TR\hat{K}^N_2(t)-\frac{\Gamma^TQ}{N}]\hat{x}_i(t)+[N\hat{K}^N_2(t)^TR\hat{K}^N_2(t)+\frac{\Gamma^TQ\Gamma}{N}]\hat{x}^{(N)}(t)+\hat{K}^N_2(t)^TR \hat{v}^N(t)    \\&+\frac{\Gamma^TQ}{N}\eta\Big\}dt+\sum_{l=1}^N\hat{z}^i_{jl}(t)dw_l(t),\\
									\hat{p}^i_j(T)=&\ -\frac{\Gamma^T_fQ_f}{N}\hat{x}_i(T)+\frac{\Gamma^T_fQ_f\Gamma_f}{N} \hat{x}^{(N)}(T)+\frac{\Gamma_f^TQ_f}{N}\eta_f,\  \forall i\in\mathcal{N},\  \forall j\in\mathcal{N}_{-i}.				   
								\end{aligned}
								\right.
							\end{split}
					\end{equation}}%
					The stationarity condition \eqref{e23101603} becomes 
					\begin{equation} \label{e23101801}
						B^T\hat{p}^i_i(t)+R[\hat{K}^N_1(t)-\hat{K}^N_2(t)]\hat{x}_i(t)+NR\hat{K}^N_2(t)\hat{x}^{(N)}(t)+R\hat{v}^N(t)=0.
					\end{equation}
					
					Next we derive the explicit solution of $(\hat{K}^N_1(\cdot),\hat{K}^N_2(\cdot),\hat{v}^N(\cdot))$ through tackling the system of FBSDEs \eqref{e23101703} by virtue of Riccati equations. Meanwhile, we will show the solvability of the Riccati equations is also sufficient for the existence of Nash equilibrium by completing the square. We first introduce the following system of Riccati equations
					\begin{subequations}\label{e101409}
						\begin{align}
							\begin{split}\label{e101409a}
								\left\{
								\begin{aligned}
									&\dot{\hat{P}}^N_1(t)=\hat{P}^N_1(t)\Upsilon \hat{P}^N_1(t)\!+\!(N\!-\!1)\hat{P}^N_2(t)\Upsilon \hat{P}^N_2(t)\!+\!(N\!-\!1)\hat{P}^N_2(t)^T\Upsilon \hat{P}^N_2(t)^T-\hat{P}^N_1(t)(A+\frac{G}{N})\\&\quad\quad\quad\quad-(A+\frac{G}{N})^T\hat{P}^N_1(t)-\frac{N-1}{N}\hat{P}^N_2(t)G-\frac{N-1}{N}G^T\hat{P}^N_2(t)^T-Q+\frac{Q\Gamma}{N}+\frac{\Gamma^TQ}{N}-\frac{\Gamma^TQ\Gamma}{N^2}, \\&\hat{P}^N_1(T)=Q_f-\frac{ \Gamma_f^TQ_f}{N}-\frac{Q_f\Gamma_f}{N}+\frac{\Gamma_f^TQ_f\Gamma_f}{N^2},
								\end{aligned}
								\right.
							\end{split}\\
							\begin{split}\label{e101409b}
								\left\{
								\begin{aligned}
									&\dot{\hat{P}}^N_2(t)=\hat{P}^N_1(t)\Upsilon \hat{P}^N_2(t)+(N-1)\hat{P}^N_2(t)^T\Upsilon \hat{P}^N_3(t)+\hat{P}^N_2(t)\Upsilon \hat{P}^N_1(t)+(N-2)\hat{P}^N_2(t)\Upsilon \hat{P}^N_2(t)-\hat{P}^N_1(t)\frac{G}{N}\\&\quad\quad\quad\quad-\frac{N-1}{N}G^T\hat{P}^N_3(t)-\hat{P}^N_2(t)(A+\frac{N-1}{N}G)-(A+\frac{G}{N})^T\hat{P}^N_2(t)-\frac{\Gamma^TQ\Gamma}{N^2}+\frac{Q\Gamma}{N} ,\quad 	\\&\hat{P}^N_2(T)=\frac{\Gamma_f^TQ_f\Gamma_f}{N^2}-\frac{Q_f\Gamma_f}{N},
								\end{aligned}
								\right.
							\end{split}\\
							\begin{split}\label{e101409c}
								\left\{
								\begin{aligned}
									&\dot{\hat{P}}^N_3(t)=\hat{P}^N_1(t)\Upsilon \hat{P}^N_3(t)\!+\!(N\!-\!2)\hat{P}^N_3(t)\Upsilon \hat{P}^N_2(t)\!+\!(N\!-\!2)\hat{P}^N_2(t)^T\Upsilon \hat{P}^N_3(t)\!+\!\hat{P}^N_3(t)\Upsilon \hat{P}^N_1(t)\!+\!\hat{P}^N_2(t)^T\Upsilon \hat{P}^N_2(t)\\&\quad\quad\quad\quad-\hat{P}^N_2(t)^T\frac{G}{N}-\frac{G^T}{N}\hat{P}^N_2(t) -\hat{P}^N_3(t)(A+\frac{N-1}{N}G) -(A+\frac{N-1}{N}G)^T\hat{P}^N_3(t) -\frac{\Gamma^TQ\Gamma}{N^2},\quad \\	&\hat{P}^N_3(T)=\frac{\Gamma_f^TQ_f\Gamma_f}{N^2},
								\end{aligned}
								\right.
							\end{split}
						\end{align}	
					\end{subequations}
					and the following system of ODEs
					\begin{subequations}\label{e101410}
						\begin{align}
							\begin{split} \label{e101410a}
								\left\{
								\begin{aligned}
									&\dot{\hat{s}}^N_1(t)=-[A+\frac{G}{N}-\Upsilon \hat{P}^N_1(t)-(N-1)\Upsilon \hat{P}^N_2(t)^T]^T\hat{s}^N_1(t)-(N-1)[\frac{G}{N}-\Upsilon \hat{P}^N_2(t)]^T\hat{s}^N_2(t)+(Q-\frac{\Gamma^TQ}{N})\eta, \\&\hat{s}^N_1(T)=(\frac{\Gamma_f^TQ_f}{N}-Q_f)\eta_f,
								\end{aligned}
								\right.
							\end{split}\\
							\begin{split} \label{e101410b}
								\left\{
								\begin{aligned}
									& \dot{\hat{s}}^N_2(t)=-[\frac{G}{N}-\Upsilon \hat{P}^N_2(t)-(N-1)\Upsilon \hat{P}^N_3(t)]^T\hat{s}^N_1(t)-[A+\frac{N-1}{N}G-\Upsilon \hat{P}^N_1(t)-(N-2)\Upsilon \hat{P}^N_2(t)]^T\hat{s}^N_2(t)\\&-\frac{\Gamma^TQ}{N}\eta,\\	&\hat{s}^N_2(T)=\frac{\Gamma_f^TQ_f}{N}\eta_f .
								\end{aligned}
								\right.
							\end{split}
						\end{align}	
					\end{subequations}
					Then the explicit solution of the centralized closed-loop Nash equilibrium  is given by the following results.
					\begin{theorem}\label{th2}
						Problem (P2) admits a centralized closed-loop Nash equilibrium $(\hat{K}_{ij}(\cdot)\in \mathcal{L}_1[0,T],\hat{v}_i(\cdot)\in \mathcal{L}_2[0,T],$ $\forall i,j\in\mathcal{N})$
						if and only if the Riccati equations \eqref{e101409} admit a solution $(\hat{P}_i(\cdot),i=1,2,3)$. 
						In this case, the centralized closed-loop Nash equilibrium $(\hat{u}_i(\cdot),\forall i\in\mathcal{N})$ of Problem (P2) is given by 
						{\setlength{\abovedisplayskip}{4pt}
							\setlength{\belowdisplayskip}{-1pt}
							\setlength{\abovedisplayshortskip}{4pt}
							\setlength{\belowdisplayshortskip}{-1pt}
							\begin{equation} \label{e122201}
								\begin{split}
									\hat{u}_i(t)=&\ \hat{K}^N_1(t)\hat{x}_i(t)+\sum_{j\in\mathcal{N}_{-i}}\hat{K}^N_2(t)\hat{x}_j(t)+\hat{v}^N(t),\quad \forall i\in\mathcal{N},
								\end{split}
						\end{equation}}
						with 
						\begin{equation}\label{e101701}
							\hat{K}^N_1(t)= \! -R^{-1} B^T\hat{P}_1^N(t),\ 
							\hat{K}^N_2(t)= \!-R^{-1} B^T\hat{P}_2^N(t),\ 
							\hat{v}^N(t)=\!	-R^{-1} B^T\hat{s}_1^N(t),\!
						\end{equation} 
						where $\hat{s}^N_1(\cdot)$ is given by \eqref{e101410a}, and the value function is given by
						\begin{equation} 
							\begin{split}
								V_i(\mathbf{x_0})\triangleq&\ {\rm J}_i(\mathbf{x}_0;\hat{K}_i(\cdot)\hat{\mathbf{x}}(\cdot)+\hat{v}_i(\cdot),\hat{K}_j(\cdot)\hat{\mathbf{x}}(\cdot)+\hat{v}_j(\cdot),\forall j\in\mathcal{N}_{-i}) \\=&\ \mathbb{E}[\langle x_{i0}, [\hat{P}^N_1(0)+(N-1)\hat{P}^N_3(0)] x_{i0}\rangle+2\langle \hat{s}^N_1(0)-\hat{s}^N_2(0),x_{i0}\rangle+2\langle N\hat{s}^N_2(0),x^{(N)}_0\rangle]+|\eta_f|^2_{Q_f}\\&+\mathbb{E}\int_0^T \Big\{\langle \sigma, [\hat{P}^N_1(t)+(N-1)\hat{P}^N_3(t)]\sigma\rangle  -\langle 2(N-1)\Upsilon \hat{s}^N_2(t)+ \Upsilon \hat{s}^N_1(t),\hat{s}^N_1(t)\rangle	+	|\eta|^2_Q \Big\}dt.
							\end{split}
						\end{equation}
					\end{theorem}
					
					\textbf{Proof} 
					See Appendix \ref{proofofth2}. $\hfill\blacksquare$

					\subsection{Discussion} In this section, for the finite-population system we have derived the feedback representation \eqref{e100804} of the centralized open-loop Nash equilibrium and the centralized closed-loop equilibrium \eqref{e122201}. It is noteworthy that the system of Riccati equations \eqref{e101201} contains four ODEs, while the system of Riccati equations  \eqref{e101409} contains three ODEs. The reason is that $(\hat{P}^N_4)^T=\hat{P}^N_2$ in the closed-loop case, as shown in the proof of Theorem \ref{th2}. Therefore the solutions to Riccati equations  \eqref{e101201} and \eqref{e101409} are significantly different.  According to Theorem \ref{closedlooppro} and Theorem \ref{th2}, we conclude that  the centralized open-loop and closed-loop Nash equilibria in the finite-population system exhibit essential distinctions.

					\section{Main Results for  Decentralized Control}\label{decentralized}
					In this section, we design  the asymptotically decentralized Nash equilibria by examining the centralized open-loop and closed-loop Nash equilibria in the infinite-population system. 
					Before presenting further results, let $P^\infty_1(\cdot)$ and $P^\infty_2(\cdot)$ be the solutions to the following Riccati equations 
					{\setlength{\abovedisplayskip}{4pt}
						\setlength{\belowdisplayskip}{5pt}
						\setlength{\abovedisplayshortskip}{4pt}
						\setlength{\belowdisplayshortskip}{5pt}
						\begin{subequations}\label{e111001}
							\begin{align}
								\begin{split} \label{e111001a}
									\left\{
									\begin{aligned}
										&\dot{P}^\infty_1(t)=P^\infty_1(t)\Upsilon P^\infty_1(t)-P^\infty_1(t)A-A^TP^\infty_1(t)-Q,\\
										&P^\infty_1(T)= Q_f,
									\end{aligned}
									\right.
								\end{split}\\
								\begin{split} \label{e111001b}
									\left\{
									\begin{aligned}
										&\dot{P}^\infty_2(t)=P^\infty_1(t)\Upsilon P^\infty_2(t)+P^\infty_2(t)\Upsilon P^\infty_1(t)+P^\infty_2(t)\Upsilon P^\infty_2(t)-A^TP^\infty_2(t)-P^\infty_2(t)(A+G)-P^\infty_1(t)G+Q\Gamma,\\
										&P^\infty_2(T)=-Q_f\Gamma_f,
									\end{aligned}
									\right.
								\end{split}
							\end{align}	
					\end{subequations}}%
					$s^\infty_1(\cdot)$  be the solutions to the following ODE 
					\begin{equation} \label{e111301}
						\begin{split} 
							\left\{
							\begin{aligned}
								&\dot{s}^\infty_1(t)=-[A^T-P^\infty_1(t)\Upsilon -P^\infty_2(t)\Upsilon ]s^\infty_1(t)+Q\eta,\quad\quad\quad\quad\quad\quad\quad\quad\quad \\&s^\infty_1(T)=-Q_f\eta_f,
							\end{aligned}
							\right.
						\end{split}
					\end{equation}
					$\check{P}^\infty_3(\cdot),\check{P}^\infty_4(\cdot)$ and $\check{s}^\infty_2(\cdot)$ be the solutions to the following ODEs
					\begin{subequations}\label{e111302}
						\begin{align}
							\begin{split}\label{e111302a} 
								\left\{
								\begin{aligned}
									&\dot{\check{P}}^\infty_3(t)=\check{P}^\infty_4(t)\Upsilon P^\infty_2(t)+\check{P}^\infty_3(t)\Upsilon P^\infty_1(t)+\check{P}^\infty_3(t)\Upsilon P^\infty_2(t)-\check{P}^\infty_4(t)G\\&\quad\quad\quad\quad-G^TP^\infty_2(t)-\check{P}^\infty_3(t)(A+G)-(A+G)^T\check{P}^\infty_3(t)-\Gamma^TQ\Gamma,\\
									&\check{P}^\infty_3(T)=\Gamma_f^TQ_f\Gamma_f,
								\end{aligned}
								\right.
							\end{split}\\
							\begin{split}\label{e111302b} 
								\left\{
								\begin{aligned}
									&\dot{\check{P}}^\infty_4(t)=\check{P}^\infty_4(t)\Upsilon P^\infty_1(t)-\check{P}^\infty_4(t)A-(A+G)^T\check{P}^\infty_4(t)-G^TP^\infty_1(t)+\Gamma^TQ\\
									&\check{P}^\infty_4(T)=-\Gamma_f^TQ_f,
								\end{aligned}
								\right.
							\end{split}\\
							\begin{split} \label{e111302c}
								\left\{
								\begin{aligned}
									&\dot{\check{s}}^\infty_2(t)=-[G^T-\check{P}^\infty_4(t)\Upsilon -\check{P}^\infty_3(t)\Upsilon ]s^\infty_1(t)-(A+G)^T\check{s}^\infty_2(t)-\Gamma^TQ\eta,\quad\\ 	&\check{s}^\infty_2(T)=\Gamma_f^TQ_f\eta_f, 
								\end{aligned}
								\right.
							\end{split}
						\end{align}	
					\end{subequations}
					and $\hat{P}^\infty_3(\cdot)$ and $\hat{s}^\infty_2(\cdot)$ be the solutions to the following ODEs
					\begin{subequations}\label{e111303}
						\begin{align}
							\begin{split} \label{e111303a}
								\left\{
								\begin{aligned}
									&\dot{\hat{P}}^\infty_3(t)=P_2^\infty(t)^T\Upsilon P_2^\infty(t)+\hat{P}_3^\infty(t)\Upsilon [P^\infty_1(t)+P^\infty_2(t)]+[P^\infty_1(t)+P^\infty_2(t)^T]\Upsilon \hat{P}_3^\infty(t)-\hat{P}^\infty_3(t)(A+G)\\&\quad\quad\quad\quad-(A+G)^T\hat{P}^\infty_3(t)-P_2^\infty(t)^TG-G^TP_2^\infty(t)-\Gamma^TQ\Gamma,\\
									&\hat{P}^\infty_3(T)=\Gamma_f^TQ_f\Gamma_f,
								\end{aligned}
								\right.
							\end{split}\\
							\begin{split} \label{e111303b}
								\left\{
								\begin{aligned}
									&\dot{\hat{s}}^\infty_2(t)=-[A+G-\Upsilon P_1^\infty(t)-\Upsilon P_2^\infty(t)]^T\hat{s}^\infty_2(t)-[G-\Upsilon P_2^\infty(t)-\Upsilon \hat{P}_3^\infty(t)]^Ts^\infty_1(t)-\Gamma^TQ\eta,\\
									&\hat{s}^\infty_2(T)=\Gamma_f^TQ_f\eta_f.
								\end{aligned}
								\right.
							\end{split}
						\end{align}	
					\end{subequations}
					Note that $Q$ is indefinite, so the Riccati equation \eqref{e111001a} may not admit a solution. If the standard Riccati equation \eqref{e111001a} and non-symmetric Riccati equation \eqref{e111001b} admit  solutions $\check{K}^\infty_2(\cdot)$, then we can uniquely solve $s^\infty_1(\cdot), \check{P}^\infty_3(\cdot),$ $\check{P}^\infty_4(\cdot),\check{s}^\infty_2(\cdot),\hat{P}_3^\infty(\cdot),\hat{s}^\infty_2(\cdot)$ from  the linear ODEs \eqref{e111302}, \eqref{e111303}.

					\subsection{Asymptotic solvability for Problem (P1)}\label{sub41}
					In this subsection, we consider the limit for the solutions of Riccati equations \eqref{e101201} and ODEs \eqref{e042401}. We assume the convex
					condition \eqref{e091402} holds and give the asymptotic solvability for Problem (P1).
					\begin{definition}
						Problem (P1) has the asymptotic solvability if there exists $\check{N}_1$ such that for all $N\geq \check{N}_1$ the system of Riccati equations \eqref{e101201} admits a solution $(\check{P}^N_i(\cdot),\forall i=1,. . .,4)$.
					\end{definition}
					However,  directly taking $N\rightarrow\infty$ is not useful because this method  will cause a loss of dynamical
					information since $(\check{P}^N_2(\cdot),\check{P}^N_3(\cdot),\check{P}^N_4(\cdot))$ can vanish when $N\rightarrow \infty$. Therefore we use the  re-scale technique proposed in \cite{b110205} to consider the limit solutions which are meaningful. Let
					\begin{equation}\label{e122402a} 
						\begin{split}
							\check{\Lambda}^N_1(\cdot)\triangleq  \check{P}^N_1(\cdot),\ \  \check{\Lambda}^N_2(\cdot)\triangleq N\check{P}^N_2(\cdot),\ \ 
							\check{\Lambda}^N_3(\cdot)\triangleq N^2\check{P}^N_3(\cdot),\ \ \check{\Lambda}^N_4(\cdot)\triangleq N\check{P}^N_4(\cdot),\ \ 
							\check{\varphi}^N_1(\cdot)\triangleq \check{s}^N_1(\cdot),\ \ 
							\check{\varphi}^N_2(\cdot)\triangleq N\check{s}^N_2(\cdot).
						\end{split}
					\end{equation}
					By some elementary estimates, we have the following results.  
					\begin{theorem}\label{theorem101901}
						Problem (P1) has the asymptotic solvability if and only if the system of Riccati equations \eqref{e111001} admits a solution $P^\infty_1(\cdot),P^\infty_2(\cdot)$. In this case, we have the following estimations.
						\begin{subequations}\label{e101706}
							\begin{align}
								\label{e101706a}\sup_{t\in[0,T]}	&(|	\check{\Lambda}^N_1(t)-P^\infty_1(t)|+|	\check{\Lambda}^N_2(t)-P^\infty_2(t)|+|\check{\varphi}^N_1(t)-s^\infty_1(t)|)=O(\frac{1}{N}),\\
								\label{e101706b}\sup_{t\in[0,T]}	&(|	\check{\Lambda}^N_3-\check{P}^\infty_3(t)|+|	\check{\Lambda}^N_4-\check{P}^\infty_4(t)|+|\check{\varphi}^N_2(t)-\check{s}^\infty_2(t)|)=O(\frac{1}{N}).
							\end{align}
						\end{subequations}
					\end{theorem}
					\textbf{Proof} 
					See Appendix \ref{proofoftheorem101901}. $\hfill\blacksquare$

					%

					\subsection{Asymptotic solvability for Problem (P2)}\label{sub42}
					Next we consider the meaning limit solutions for Riccati equations \eqref{e101409} and ODEs \eqref{e101410} by the re-scale technique. The asymptotic solvability for Problem (P2) is given as follows. 
					\begin{definition}
						Problem (P2) has the asymptotic solvability if there exists $\hat{N}_1$ such that for all $N\geq \hat{N}_1$ the system of Riccati equations \eqref{e101409} admits a solution $(\hat{P}^N_i(\cdot),\forall i=1,. . .,3)$.
					\end{definition}
					Let 
					\begin{equation} \label{e122402b}
						\begin{split}
							\hat{\Lambda}^N_1(\cdot)\triangleq  \hat{P}^N_1(\cdot),\quad \hat{\Lambda}^N_2(\cdot)\triangleq N\hat{P}^N_2(\cdot),\quad 
							\hat{\Lambda}^N_3(\cdot)\triangleq N^2\hat{P}^N_3(\cdot),\quad 
							\hat{\varphi}^N_1(\cdot)\triangleq \hat{s}^N_1(\cdot),\quad
							\hat{\varphi}^N_2(\cdot)\triangleq N\hat{s}^N_2(\cdot).
						\end{split}
					\end{equation}
					Similar to Theorem \ref{theorem101901}, we have the following results.
					\begin{theorem}\label{theorem101902}
						Problem (P2) has the asymptotic solvability if and only if the system of Riccati equations \eqref{e111001} admits a solution $P^\infty_1(\cdot),P^\infty_2(\cdot)$. In this case, we have the following estimations. 
						\begin{subequations}\label{e101803}
							\begin{align}
								\label{e101803a}		\sup_{t\in[0,T]}	&(|\hat{\Lambda}^N_1(t)-P_1^\infty(t)|+|\hat{\Lambda}^N_2(t)-P_2^\infty(t)|+|\hat{\varphi}^N_1-s_1^\infty(t)|)=O(\frac{1}{N}), \\
								\label{e101803b}		\sup_{t\in[0,T]}	&(|\hat{\Lambda}^N_3(t)-\hat{P}_3^\infty(t)|+|\hat{\varphi}^N_2(t)-\hat{s}_2^\infty(t)|)=O(\frac{1}{N}).
							\end{align}
						\end{subequations}
					\end{theorem}
					\textbf{Proof} 
					See Appendix \ref{proofoftheorem101902}.  $\hfill\blacksquare$

					\subsection{Asymptotically decentralized Nash equilibria for Problem (P3) and Problem (P4)}  By the re-scale technique, we have studied the meaningful limit solutions to the Riccati equations \eqref{e101201}, \eqref{e101409} and ODEs \eqref{e042401}, \eqref{e101410} in Theorem \ref{theorem101901} and Theorem \ref{theorem101902}. The results show that the limit solutions of $(\check{\Lambda}^N_1(\cdot),\check{\Lambda}^N_2(\cdot),\check{\varphi}^N_1(\cdot))$ in \eqref{e122402a} and $(\hat{\Lambda}^N_1(\cdot),\hat{\Lambda}^N_2(\cdot),\hat{\varphi}^N_1(\cdot))$ in \eqref{e122402b} are both $(P^\infty_1(\cdot),P^\infty_2(\cdot),s^\infty_1(\cdot))$. In this subsection, we use the mean field approximation to design the asymptotically decentralized Nash equilibria based on the centralized open-loop and closed-loop Nash equilibria.
					
					First, we solve the dynamics of the mean field term under the centralized open-loop and closed-loop Nash equilibria, respectively. The feedback representation \eqref{e100804} can be rewritten as 
					\begin{equation} \label{e122301}
						\begin{split}
							\check{u}_i(t)=&\ [\check{K}^N_1(t)-\check{K}^N_2(t)]\check{x}_i(t)+N\check{K}^N_2(t)\check{x}^{(N)}(t)+\check{v}^N(t),\quad \forall i\in\mathcal{N}.
						\end{split}
					\end{equation}
					Under the feedback representation \eqref{e122301}, the closed-loop dynamics of $\check{x}_i(\cdot)$, $\forall i\in\mathcal{N}$, is given by
					\begin{equation}\label{e111601}
						\begin{split}
							d\check{x}_i(t) =\Big\{[A+B\check{K}^N_1(t)-B\check{K}^N_2(t)]\check{x}_i(t) +[G+NB\check{K}^N_2(t)]\check{x}^{(N)}(t) +B\check{v}^N(t)  \Big\}dt+\sigma dw_i(t) ,\quad \check{x}_i(0)=x_{i0}.
						\end{split}
					\end{equation}
					Therefore the dynamics of the mean field term $\check{x}^{(N)}(\cdot)$ under the feedback representation \eqref{e122301} is given by
					\begin{equation}\label{e101903}
						\begin{split}
							d\check{x}^{(N)}(t)= \Big\{[A+G+B\check{K}^N_1(t)+(N-1)B\check{K}^N_2(t)]\check{x}^{(N)}(t) +B\check{v}^N(t) \Big\}dt+\sum_{i=1}^N\frac{\sigma}{N} dw_i(t),\quad \check{x}^{(N)}(0)=x^{(N)}_{0}.
						\end{split}
					\end{equation}
					Similar to \eqref{e122301}, \eqref{e111601} and \eqref{e101903},  the closed-loop dynamics of the mean field term $\hat{x}^{(N)}(\cdot)$ under the centralized closed-loop Nash equilibrium \eqref{e122201} is given by 
					\begin{align}
						\begin{split}\label{e101906}
							d\hat{x}^{(N)}(t)= \Big\{[A+G+B\hat{K}^N_1(t)+(N-1)B\hat{K}^N_2(t)]\hat{x}^{(N)}(t) +B \hat{v}^N(t) \Big\}dt+\sum_{i=1}^N\frac{\sigma}{N} dw_i(t),\quad \hat{x}^{(N)}(0)=x^{(N)}_{0}.
						\end{split}
					\end{align}
					Next we consider the limit for  $\check{x}^{(N)}(\cdot)$ and $\hat{x}^{(N)}(\cdot)$ by letting $N\rightarrow\infty$. Let 
					\begin{equation}\label{e122405}
						K^\infty_1(\cdot)=-BR^{-1}P^\infty_1(\cdot),\quad K^\infty_2(\cdot)=-BR^{-1}P^\infty_2(\cdot),\quad 
						\varphi^\infty_1(\cdot)=-BR^{-1}s^\infty_1(\cdot).
					\end{equation}
					Then we have the following results.
					\begin{proposition}\label{t122401}
						Suppose Problem (P1) and Problem (P2) have the asymptotic solvability. Then we have the following estimations. 
						\begin{subequations}\label{e122404}
							\begin{align}
								\sup_{t\in[0,T]}	&(|\check{K}^N_1(t)-K^\infty_1(t)|+|N\check{K}^N_2(t)-K^\infty_2(t)|+|\check{\varphi}^N_1-\varphi^\infty_1(t)|)=O(\frac{1}{N}), \\
								\sup_{t\in[0,T]}	&(|\hat{K}^N_1(t)-K^\infty_1(t)|+|N\hat{K}^N_2(t)-K^\infty_2(t)|+|\hat{\varphi}^N_1-\varphi^\infty_1(t)|)=O(\frac{1}{N}). 
							\end{align}
						\end{subequations}
						
					\end{proposition}
					\textbf{Proof} 
					The proposition follows by  the definitions of $(\check{K}^N_1,\check{K}^N_2,\check{s}^N_1)$ in \eqref{e122403a}, $(\hat{K}^N_1,$ $\hat{K}^N_2,\hat{s}^N_1)$ in \eqref{e101701}, $(K^\infty_1,K^\infty_2,\varphi^\infty_1)$ in \eqref{e122404}, and Theorems \ref{theorem101901} and \ref{theorem101902}. $\hfill\blacksquare$ 
					
					To denote the limit of \eqref{e101903} and \eqref{e101906}  when $N\rightarrow \infty$, we introduce the  following mean field dynamics
					\begin{equation}\label{e101902}
						\begin{split}
							\dot{\bar{x}}(t)=[A+G+BK^\infty_1(t)+BK^\infty_2(t)]\bar{x}(t)+B\varphi_1^\infty(t) ,\quad \bar{x}(0)=\bar{x}_0.
						\end{split}
					\end{equation}
					For the mean field approximation, we have the following results. 
					\begin{proposition}\label{t122301}
						Let (A1) holds. Suppose Problem (P1) and Problem (P2) have the asymptotic solvability. Then we have
						\begin{subequations}
							\begin{align}
								\label{e101706c}	\sup_{t\in[0,T]}&\mathbb{E}|\check{x}^{(N)}(t)-\bar{x}(t)|^2=O(\frac{1}{N}),\\
								\label{e101803c}		\sup_{t\in[0,T]}&\mathbb{E}|\hat{x}^{(N)}(t)-\bar{x}(t)|^2=O(\frac{1}{N}). 
							\end{align}
						\end{subequations}
					\end{proposition}
					\textbf{Proof} 
					See Appendix \ref{proofoft122301}. $\hfill\blacksquare$

					Based on Propositions \ref{t122401} and \ref{t122301}, from \eqref{e100804} and \eqref{e122201}, we design the asymptotically decentralized Nash equilibrium 
					\begin{equation}\label{e122406} 
						\begin{split}
							u_i^d(t)= K^\infty_1(t)x_i^d(t)+K^\infty_2(t)\bar{x}(t)+\varphi^\infty_1(t),\quad i\in\mathcal{N},
						\end{split}
					\end{equation}
					where $x^d_i(\cdot)$ is the realized state under the decentralized control $u_i^d(\cdot)$, $K^\infty_1(\cdot),K^\infty_2(\cdot),$ $\varphi^\infty_1(\cdot)$ and $\bar{x}(\cdot)$ are given by \eqref{e122405} and \eqref{e101902}, respectively. 
					\subsection{Discussion} In this section, the meaningful limit solutions to the Riccati equations \eqref{e101201}, \eqref{e101409} and ODEs \eqref{e042401}, \eqref{e101410} have been studied in Theorem \ref{theorem101901} and Theorem \ref{theorem101902} by utilizing the re-scale technique. By examining the feedback representation \eqref{e100804} of the centralized closed-loop Nash equilibrium and the centralized closed-loop Nash equilibrium \eqref{e122201} in the infinite-population system, we obtain the same asymptotically decentralized Nash equilibrium  \eqref{e122406}.  We can compute $P^\infty_1(\cdot),P^\infty_2(\cdot),s^\infty_1(\cdot)$ and $\bar{x}(\cdot)$ by using \eqref{e111001}, \eqref{e111301} and \eqref{e101902} offline.  Therefore, each player only needs to use their own state information to implement the decentralized control \eqref{e122406}. Therefore, while open-loop and closed-loop Nash equilibria are significantly different in the finite-population system, they   asymptotically coincide in the infinite-population system. The results indicate that the size of the population has a essential role on the open-loop and closed-loop Nash equilibria, while the choice of open-loop or closed-loop strategies does not influence the design of asymptotically decentralized control in LQ MFGs.

					\section{Conclusion}\label{conclusion}
					
					In this paper, we adopt a unified framework to concurrently study the (asymptotically) open-loop and closed-loop Nash equilibria for LQ MFGs in both finite and infinite-population systems. Our aim is to investigate the impact of population scale on the design of open-loop and closed-loop strategies and to explore whether adopting different concepts of open-loop and closed-loop strategies will affect the design of asymptotically decentralized Nash equilibria. The results show that centralized open-loop and closed-loop Nash equilibria are significantly different in the finite-population system. Via the re-scale technique we study the meaningful limit solutions of the system of Riccati equations in the infinite-population system. The asymptotically decentralized Nash equilibria that only utilize the players' own state information are designed by examining centralized open-loop and closed-loop Nash equilibria in the infinite-population system. The results demonstrate that the designed decentralized control based on centralized open-loop and closed-loop Nash equilibria are consistent. Therefore, the choice between open-loop and closed-loop strategies bears no consequence on the design of decentralized control in LQ MFGs.

					\appendix

					\section{Proofs for Results in Section \ref{centralized}}

					
					\subsection{Proof of Theorem \ref{th1}}\label{proofofthe1} Henceforward, for notation simplicity, we will drop time variable $t$ if no confusion in the proof. By Definition \ref{def1}, $(\check{u}_i(\cdot)\in \mathcal{U}_c[0,T],\forall i\in\mathcal{N})$ 
					is an open-loop Nash equilibrium if and only if the following conditions hold  $\forall\theta  \in \mathbb{R}$:
					\begin{equation}\label{e101301}
						\begin{split}
							{\rm J}_i(\mathbf{x}_0;\check{u}_i(\cdot)+\theta  u_i(\cdot),\check{u}_{-i}(\cdot))\geq 	{\rm J}_i(\mathbf{x}_0;\check{u}_i(\cdot),\check{u}_{-i}(\cdot))  ,\quad \forall i\in\mathcal{N}.
						\end{split}
					\end{equation}
					Under the control $(\check{u}_i(\cdot),\forall i\in\mathcal{N})$, the state $\check{x}_i(\cdot)$ satisfies
					\begin{equation}\label{e101303}
						d\check{x}_i =(A\check{x}_i +G\check{x}^{(N)} +B\check{u}_i )dt+\sigma dw_i ,\quad \check{x}_i(0)=x_{i0},\quad \forall i\in\mathcal{N}.
					\end{equation}	
					For any $i\in\mathcal{N}$, $u_i(\cdot)\in \mathcal{U}_c[0,T]$ and $\theta \in \mathbb{R}$, let $(x_j^{ i,\theta}(\cdot), \forall j\in\mathcal{N})$ be the solution of the following system of perturbed state equations:
					\begin{subequations}\label{e101302}
						\begin{align}
							dx_i^{i,\theta} =&\ (Ax_i^{i,\theta} +Gx^{i,(N),\theta} +B\check{u}_i +\theta Bu_i )dt+\sigma dw_i ,\quad x_i^{i,\theta}(0)=x_{i0} \\
							dx_j^{i,\theta} =&\ (Ax_j^{i,\theta} +Gx^{i,(N),\theta} +B\check{u}_j )dt+\sigma dw_j ,\quad x_i^{j,\theta}(0)=x_{j0},\quad \forall  j\in\mathcal{N}_{-i},
						\end{align}
					\end{subequations}
					where $x^{i,(N),\theta}(\cdot)\triangleq (1/N)\sum_{j=1}^Nx_j^{i,\theta}(\cdot)$. Then by \eqref{e101303} and \eqref{e101302}, $\check{y}_j^i(\cdot)\triangleq(x_j^{i,\theta}(\cdot)-\check{x}_{j}(\cdot))/\theta,\forall j\in\mathcal{N}$, is independent of $\theta$ satisfying \eqref{e101305}. Let $(\check{p}_i^j(\cdot),\check{z}^i_{jk}(\cdot),\forall i,j,k\in\mathcal{N})$ be the solution of the system of backward stochastic differential equations in \eqref{e101101}. 
					By \eqref{e101101}, \eqref{e101305} and It\^{o}' formula to $t\mapsto \mathbb{E}\langle\check{p}_j^i,\check{y}_j^i\rangle,\forall j\in\mathcal{N}$, we have
					\begin{subequations}
						\begin{align}
							\begin{split}
								\mathbb{E}	\langle \check{p}_{i}^i(T),\check{y}_i^i(T) \rangle=\	&\mathbb{E}\int_0^T\Big\{\langle  -[G^T\check{p}^{i,(N)} +(Q-\frac{\Gamma^TQ}{N})\check{x}_i -(Q\Gamma-\frac{\Gamma^TQ\Gamma}{N}) \check{x}^{(N)}-(Q-\frac{\Gamma^TQ}{N})\eta], \check{y}_i^i \rangle\\&+\langle G^T\check{p}_i^i , \check{y}^{i,(N)} \rangle +\langle B^T\check{p}_i^i , u_i \rangle\Big\},
							\end{split}\\
							\begin{split}
								\mathbb{E}\langle \check{p}_{j}^i(T),\check{y}_j^i(T) \rangle
								=\ &\mathbb{E}\int_0^T\Big[\langle -(G^T\check{p}^{i,(N)} -\frac{\Gamma^TQ}{N}\check{x}_i +\frac{\Gamma^TQ\Gamma}{N}\check{x}^{(N)} +\frac{\Gamma^TQ}{N}\eta), \check{y}_j^i \rangle+\langle G^T\check{p}^i_j ,\check{y}^{i,(N)} \rangle \Big],\quad j\in\mathcal{N}_{-i},
							\end{split}
						\end{align}
					\end{subequations}
					which further gives
					\begin{align}
						\begin{split}\label{e101306}
							\sum_{j=1}^N \mathbb{E}	\langle \check{p}_{j}^i(T),\check{y}_j^i(T) \rangle
							=&\ \mathbb{E}\int_0^T\Big\{\langle B^T\check{p}_i^i , u_i \rangle-\sum_{j\in\mathcal{N}_{-i}}\langle  [-\frac{\Gamma^TQ}{N}\check{x}_i +\frac{\Gamma^TQ\Gamma}{N}\check{x}^{(N)} +\frac{\Gamma^TQ}{N}\eta], \check{y}_j^i \rangle\\&-\langle  (Q-\frac{\Gamma^TQ}{N})\check{x}_i -(Q\Gamma-\frac{\Gamma^TQ\Gamma}{N}) \check{x}^{(N)} -(Q-\frac{\Gamma^TQ}{N})\eta, \check{y}_i^i \rangle \Big\}.
						\end{split}
					\end{align}
					The variation of the cost functional for $\mathcal{A}_i,i\in\mathcal{N},$  is 
					\begin{align}
						\nonumber	&{\rm J}_i(\mathbf{x}_0;\check{u}_i(\cdot)+\theta u_i(\cdot),\check{u}_j(\cdot),\forall j\in\mathcal{N}_{-i})-{\rm J}_i(\mathbf{x}_0;\check{u}_i(\cdot),\check{u}_j(\cdot),\forall j\in\mathcal{N}_{-i})\\
						\label{e111304}=&\ \theta^2\Big[ \mathbb{E}\int_0^T\big(|\check{y}_i^i -\Gamma \check{y}^{i,(N)} |^2_Q+|\check{u}_i |^2_R\big)dt+\mathbb{E}|\check{y}_i^i(T)-\Gamma_f \check{y}^{i,(N)}(T)|^2_{Q_f} \Big]\\&+2\theta \Big[ \mathbb{E}\int_0^T \big(\langle \check{x}_i -\Gamma \check{x}^{(N)} -\eta, Q\check{y}_i^i -Q\Gamma \check{y}^{i,(N)} \rangle+\langle R\check{u}_i ,  u_i \rangle\big)dt\nonumber \\&+\mathbb{E} \langle Q_f\check{x}_i(T)-Q_f\Gamma_f \check{x}^{(N)}(T)-Q_f\eta_f, \check{y}_i^i(T)-\Gamma_f \check{y}^{i,(N)}(T)\rangle  \Big].\nonumber
					\end{align}
					Adding \eqref{e101306} to \eqref{e111304}, we obtain 
					\begin{equation}
						\begin{split}
							&{\rm J}_i(\mathbf{x}_0;\check{u}_i(\cdot)+\theta u_i(\cdot),\check{u}_j(\cdot),\forall j\in\mathcal{N}_{-i})-{\rm J}_i(\mathbf{x}_0;\check{u}_i(\cdot),\check{u}_j(\cdot),\forall j\in\mathcal{N}_{-i}) \\
							=&\ \theta^2\Big\{ \mathbb{E}\int_0^T\big(|\check{y}_i^i -\Gamma \check{y}^{i,(N)} |^2_Q+|\check{u}_i |^2_R\big)dt+\mathbb{E}|\check{y}_i^i(T)-\Gamma_f \check{y}^{i,(N)}(T)|^2_{Q_f} \Big\}+2\theta \mathbb{E}\int_0^T\langle R\check{u}_i +B^T\check{p}_i^i ,  u_i \rangle dt. 
						\end{split}
					\end{equation} 
					Therefore \eqref{e101301} holds if and only if the convexity condition \eqref{e091402} holds and $\check{u}_i=-R^{-1} B^T\check{p}_i^i$. The proof is completed.

					\subsection{Proof of Theorem \ref{closedlooppro}}\label{proofofclosedlooppro}
					Let
					\begin{subequations}\label{e101103}
						\begin{align}
							\check{\chi}^N_1 \triangleq&\ \check{p}^i_i -(\check{P}^N_1 -\check{P}^N_2 )\check{x}_i -N\check{P}^N_2 \check{x}^{(N)} ,\\
							\check{\chi}^N_2 \triangleq &\ \check{p}^i_j -(\check{P}^N_4 -\check{P}^N_3 )\check{x}_i -N\check{P}^N_3 \check{x}^{(N)} ,
						\end{align}
					\end{subequations}
					where $\check{P}^N_i(\cdot):[0,T]\rightarrow\mathbb{R}^{n\times n},i=1, . . .,4,$ is the solution of the Riccati equations \eqref{e101201}, and $(\check{x}_i(\cdot),\check{p}_i^j(\cdot),\check{z}_{jk}^i(\cdot),\forall i,j,k\in\mathcal{N})$ is the adapted solution of the system of FBSDEs \eqref{e101101}. To match the terminal condition, we impose the requirements
					\begin{equation} 
						\begin{split}
							\check{\chi}^N_1(T)=-(Q_f-\frac{\Gamma^T_fQ_f}{N})\eta_f,\quad 
							\check{\chi}^N_2(T)= \frac{\Gamma^T_fQ_f}{N}\eta_f .
						\end{split}
					\end{equation}
					By \eqref{e101102}, we have
					\begin{equation}\label{e100801}
						\begin{split}
							\check{p}^{i,(N)} =  \Big[\frac{1}{N}(\check{P}^N_1 -\check{P}^N_2)+\frac{N-1}{N}(\check{P}^N_3 -\check{P}^N_4)  \Big]\check{x}_i  +[\check{P}^N_2 +(N-1)\check{P}^N_3 ]\check{x}^{(N)} +\frac{1}{N}\check{\chi}_1^N+\frac{N-1}{N}\check{\chi}_2^N.
						\end{split}
					\end{equation}
					It follows by \eqref{e101102} and \eqref{e101103} that
					\begin{equation}\label{e101104}
						\begin{split}
							\check{u}_i =-R^{-1}B^T(\check{P}^N_1 -\check{P}^N_2 )\check{x}_i -NR^{-1}B^T\check{P}^N_2 \check{x}^{(N)} -R^{-1}B^T\check{\chi}^N_1 ,\quad \forall i\in\mathcal{N}.
						\end{split}
					\end{equation}
					Therefore $\check{x}_i(\cdot)$, $\forall i\in\mathcal{N}$, satisfies
					\begin{equation}\label{e100802}
						\begin{split}
							d\check{x}_i =[A\check{x}_i -\Upsilon (\check{P}^N_1 -\check{P}^N_2 )\check{x}_i -\Upsilon \check{\chi}_1 +(G-N\Upsilon \check{P}^N_2 )\check{x}^{(N)} ]dt+\sigma dw_i , 
						\end{split}
					\end{equation}	
					with $\check{x}_i(0)=x_{i0}$,	which further gives
					\begin{equation}\label{e101901}
						\begin{split}
							d\check{x}^{(N)} =&\ \Big\{(A+G)\check{x}^{(N)} -\Upsilon [\check{P}^N_1 +(N-1)\check{P}^N_2 ]\check{x}^{(N)} -\Upsilon \check{\chi}_1 \Big\}dt+\sum_{i=1}^N\frac{\sigma}{N} dw_i , 
						\end{split}
					\end{equation}
					with $\check{x}^{(N)}(0)=x^{(N)}_0$. By \eqref{e101101}, \eqref{e101201}, \eqref{e100801}, \eqref{e101901}, \eqref{e100802}  and applying It\^{o}'s formula to $	t\mapsto \check{p}^i_i -(\check{P}^N_1 -\check{P}^N_2)\check{x}_i -N\check{P}^N_2 \check{x}^{(N)}$ and  $t\mapsto \check{p}^i_j -(\check{P}^N_4 -\check{P}^N_3 )\check{x}_i -N\check{P}^N_3 \check{x}^{(N)}$,
					we obtain
					\begin{equation}\label{e100803}
						\check{z}^i_{ii} =\check{P}^N_1 \sigma,\quad  	\check{z}^i_{ik} =\check{P}^N_2 \sigma,\quad 
						\check{z}^i_{ji} =\check{P}^N_4 \sigma,\quad 
						\check{z}^i_{jk} =\check{P}^N_3 \sigma, \ \ \forall i\in\mathcal{N},\ \  \forall j,k\in\mathcal{N}_{-i},
					\end{equation}
					and 
					\begin{subequations} 
						\begin{align}
							\begin{split}
								\left\{
								\begin{aligned}
									&\dot{\check{\chi}}^N_1 =-[A+\frac{G}{N}-\Upsilon(\check{P}^N_1)^T-(N-1)\Upsilon(\check{P}^N_2)^T]^T\check{\chi}^N_1 -\frac{N-1}{N}G^T\check{\chi}^N_2 +(Q-\frac{\Gamma^TQ}{N})\eta,\quad \\&\check{\chi}^N_1(T)=(\frac{\Gamma_f^TQ_f}{N}-Q_f)\eta_f,
								\end{aligned}
								\right.
							\end{split}\\
							\begin{split}
								\left\{
								\begin{aligned}
									&\dot{\check{\chi}}^N_2 =-[\frac{G}{N}-\Upsilon(\check{P}^N_4)^T-(N-1)\Upsilon(\check{P}^N_3)^T]^T\check{\chi}^N_1 -(A+\frac{N-1}{N}G)^T\check{\chi}^N_2  -\frac{\Gamma^TQ}{N}\eta,\quad\\ 	&\check{\chi}^N_2(T)=\frac{\Gamma_f^TQ_f}{N}\eta_f ,
								\end{aligned}
								\right.
							\end{split}
						\end{align}	 
					\end{subequations}
					which gives $\check{\chi}^N_i(\cdot)=\check{s}^N_i(\cdot), i=1,2.$ Using $\check{s}_1(\cdot)$ to replace $\check{\chi}_1(\cdot )$ in \eqref{e101104}, we obtain \eqref{e100804}. The proof is completed. 
					

					\subsection{Proof of Theorem \ref{th2}}\label{proofofth2}
					We first prove the necessity. Assume  Problem (P2) admits a centralized closed-loop Nash equilibrium  $(\hat{K}_{ij}(\cdot),\hat{v}_i(\cdot),$ $\forall i,j\in\mathcal{N})$. We will show that the Riccati equations \eqref{e101409} admits solutions. Under the given parameters $(\hat{K}^N_1(\cdot),\hat{K}^N_2(\cdot))$, the following system of ODEs  have a solution $(\hat{P}^N_i , i=1,2,3,4)$ on $t \in[0,T]$:
					\begin{subequations}\label{e101601}
						\begin{align}
							\begin{split} 
								\left\{
								\begin{aligned}
									&\dot{\hat{P}}^N_1 =-\hat{P}^N_1 B\hat{K}^N_1 -(\hat{K}^N_1)^TB^T\hat{P}^N_1 -(\hat{K}^N_1)^TR\hat{K}^N_1 -(N-1)\hat{P}^N_2 B\hat{K}^N_2 -(N-1)(\hat{K}^N_2 )^TB^T\hat{P}^N_4 \\&\quad\quad\quad-\hat{P}^N_1 (A+\frac{G}{N})-(A+\frac{G}{N})^T\hat{P}^N_1 -\frac{N-1}{N}\hat{P}^N_2 G-\frac{N-1}{N}G^T\hat{P}^N_4 -Q+\frac{Q\Gamma}{N}+\frac{\Gamma^TQ}{N}-\frac{\Gamma^TQ\Gamma}{N^2}, \\&\hat{P}^N_1(T) =Q_f-\frac{ \Gamma_f^TQ_f}{N}-\frac{Q_f\Gamma_f}{N}+\frac{\Gamma_f^TQ_f\Gamma_f}{N^2},
								\end{aligned}
								\right.
							\end{split}\\
							\begin{split} 
								\left\{
								\begin{aligned}
									&\dot{\hat{P}}^N_2 =-\hat{P}^N_1 B\hat{K}^N_2 -\hat{P}^N_2 B\hat{K}^N_1 -(\hat{K}^N_1)^TB^T\hat{P}^N_2 -(\hat{K}^N_1)^TR\hat{K}^N_2 -(N-1)(\hat{K}^N_2 )^TB^T\hat{P}^N_3-\hat{P}^N_1\frac{G}{N}\\&\quad\quad\quad -(N-2)\hat{P}^N_2 B\hat{K}^N_2 -\hat{P}^N_2 (A+\frac{N-1}{N}G) -(A+\frac{G}{N})^T\hat{P}^N_2 -\frac{N-1}{N}G^T\hat{P}^N_3 -\frac{\Gamma^TQ\Gamma}{N^2}+\frac{Q\Gamma}{N} ,\quad 	\\&\hat{P}^N_2(T) =\frac{\Gamma_f^TQ_f\Gamma_f}{N^2}-\frac{Q_f\Gamma_f}{N},
								\end{aligned}
								\right.
							\end{split}\\
							\begin{split} 
								\left\{
								\begin{aligned}
									&\dot{\hat{P}}^N_3 =-(\hat{K}^N_2 )^TR\hat{K}^N_2 -(\hat{K}^N_1)^TB^T\hat{P}^N_3 -\hat{P}^N_3 B\hat{K}^N_1 -\hat{P}^N_4 B\hat{K}^N_2-\hat{P}^N_4 \frac{G}{N}  -(\hat{K}^N_2 )^TB^T\hat{P}^N_2-\frac{G^T}{N}\hat{P}^N_2\\&\quad\quad\quad-(N\!-\!2)\hat{P}^N_3 B\hat{K}^N_2 \!-\!(N\!-\!2)(\hat{K}^N_2 )^TB^T\hat{P}^N_3 \! -\!\hat{P}^N_3 (A\!+\!\frac{N-1}{N}G) \!-\!(A\!+\!\frac{N-1}{N}G)^T\hat{P}^N_3  -\frac{\Gamma^TQ\Gamma}{N^2},\quad \\	&\hat{P}^N_3(T) =\frac{\Gamma_f^TQ_f\Gamma_f}{N^2},
								\end{aligned}
								\right.
							\end{split}\\
							\begin{split} 
								\left\{
								\begin{aligned}
									&\dot{\hat{P}}^N_4 =-(\hat{K}^N_2 )^TR\hat{K}^N_1 -\hat{P}^N_4 B\hat{K}^N_1 -\hat{P}^N_3 B\hat{K}^N_2 -(\hat{K}^N_1)^TB^T\hat{P}^N_4 -\frac{G^T}{N}\hat{P}^N_1-(\hat{K}^N_2 )^TB^T\hat{P}^N_1-	\hat{P}^N_4 (A\!+\!\frac{G}{N})\! \\&\quad\quad\quad -(N\!-\!2)(\hat{K}^N_2 )^TB^T\hat{P}^N_4 -(N\!-\!2)\hat{P}^N_3 B\hat{K}^N_2  -\!\frac{N-1}{N}\hat{P}^N_3 G \!-\!(A\!+\!\frac{N-1}{N}G)^T\hat{P}^N_4 \!+\!\frac{\Gamma^TQ}{N} \!-\!\frac{\Gamma^TQ\Gamma}{N^2}, 	\\&\hat{P}^N_4(T) =\frac{\Gamma_f^TQ_f\Gamma_f}{N^2}-\frac{\Gamma_f^TQ_f}{N},
								\end{aligned}
								\right.
							\end{split}
						\end{align}	
					\end{subequations}
					and the following system of ODEs have solution $(\hat{s}^N_i ,i=1,2)$ on $t\in[0,T]$:
					\begin{subequations}\label{e101602}
						\begin{align}
							\begin{split}
								\left\{
								\begin{aligned}
									&\dot{\hat{s}}^N_1 =-(A+\frac{G}{N}
									+B\hat{K}^N_1)^T\hat{s}^N_1-(N-1) (\frac{G}{N}+B\hat{K}^N_2 )^T\hat{s}^N_2  \\&\quad\quad\quad-[\hat{P}^N_1 B+(N-1)\hat{P}^N_2 B+(\hat{K}^N_1)^TR]\hat{v}^N  +(Q-\frac{\Gamma^TQ}{N})\eta, \\&\hat{s}^N_1(T) =(\frac{\Gamma_f^TQ_f}{N}-Q_f)\eta_f,
								\end{aligned}
								\right.
							\end{split}\\
							\begin{split}
								\left\{
								\begin{aligned}
									& \dot{\hat{s}}^N_2 =-(\frac{G}{N}+B\hat{K}^N_2 )^T\hat{s}^N_1  -[A+\frac{N-1}{N}G+B\hat{K}^N_1 -(N-1)B\hat{K}^N_2 ]^T\hat{s}^N_2 \\&\quad\quad\quad-[\hat{P}^N_4B +(N-1)\hat{P}^N_3B+(\hat{K}^N_2 )^TR]\hat{v}^N    -\frac{\Gamma^TQ}{N}\eta,\\	&\hat{s}^N_2(T)=\frac{\Gamma_f^TQ_f}{N}\eta_f .
								\end{aligned}
								\right.
							\end{split}
						\end{align}	
					\end{subequations}
					By Proposition \ref{prop2} we have the system of  FBSDEs \eqref{e23101703} and the stationarity condition \eqref{e23101801}.  It follows by \eqref{e23101703} that
					\begin{equation}\label{e101603}
						\begin{split}
							d\hat{x}^{(N)}=
							[(A+G+B\hat{K}^N_1 +(N-1)B\hat{K}^N_2 )\hat{x}^{(N)} +B\hat{v}^N ]dt+\sum_{i=1}^N\frac{\sigma}{N} dw_i ,
						\end{split}
					\end{equation}
					with $\hat{x}^{(N)}(0)=x^{(N)}_0$.  By \eqref{e23101703}, \eqref{e101601}, \eqref{e101602}, \eqref{e101603} and applying It\^{o}'s formula to $	t\mapsto (\hat{P}^N_1 -\hat{P}^N_2)\hat{x}_i +N\hat{P}^N_2 \hat{x}^{(N)} +\hat{s}^N_1$ and $t\mapsto (\hat{P}^N_4 -\hat{P}^N_3 )\hat{x}_i +N\hat{P}^N_3 \hat{x}^{(N)} +\hat{s}^N_2$,
					we obtain
					\begin{subequations}\label{e23110701}
						\begin{align}
							&d\big((\hat{P}^N_1 -\hat{P}^N_2 )\hat{x}_i +N\hat{P}^N_2 \hat{x}^{(N)} +\hat{s}^N_1 \big)	\nonumber\\
							=&\ (\dot{\hat{P}}^N_1 -\dot{\hat{P}}^N_2 )\hat{x}_i dt+N\dot{\hat{P}}^N_2 \hat{x}^{(N)} dt+\dot{\hat{s}}^N_1 dt+(\hat{P}^N_1 -\hat{P}^N_2 )d\hat{x}_i +N\hat{P}^N_2 d\hat{x}^{(N)} \nonumber\\
							\label{e23101901}=&-\Big\{(A+B\hat{K}^N_1 -B\hat{K}^N_2 )^T[(\hat{P}^N_1 -\hat{P}^N_2 )\hat{x}_i +N\hat{P}^N_2 \hat{x}^{(N)} +\hat{s}^N_1 ]+(G+NB\hat{K}^N_2 )^T\Big[[\frac{1}{N}(\hat{P}^N_1 -\hat{P}^N_2)\\&+\frac{N-1}{N}(\hat{P}^N_4 -\hat{P}_3) ]\hat{x}_i +[\hat{P}^N_2 +(N-1)\hat{P}^N_3]\hat{x}^{(N)} +\frac{1}{N}\hat{s}^N_1+\frac{N-1}{N}\hat{s}^N_2 \Big]  +[(\hat{K}^N_1)^TR(\hat{K}^N_1-\hat{K}^N_2) +Q-\frac{\Gamma^TQ}{N}]	\hat{x}_i \nonumber \\&+[N(\hat{K}^N_1)^TR\hat{K}^N_2 -Q\Gamma+\frac{\Gamma^TQ\Gamma}{N}]\hat{x}^{(N)} +(\hat{K}^N_1)^TR \hat{v}^N   -(Q-\frac{\Gamma^TQ}{N})\eta\Big\}dt+(\hat{P}^N_1 -\hat{P}^N_2 )\sigma dw_i +\sum_{i=1}^N\hat{P}^N_2 \sigma dw_i ,\nonumber	  \\
							&d\big((\hat{P}^N_4 -\hat{P}^N_3 )\hat{x}_i +N\hat{P}^N_3 \hat{x}^{(N)} +\hat{s}^N_2 \big)\nonumber \\
							=&\ (\dot{\hat{P}}^N_4 -\dot{\hat{P}}^N_3 )\hat{x}_i dt+N\dot{\hat{P}}^N_3 \hat{x}^{(N)} dt+\dot{\hat{s}}^N_2 dt+(\hat{P}^N_4 -\hat{P}^N_3)d\hat{x}_i +N\hat{P}^N_3 d\hat{x}^{(N)} \nonumber\\
							\label{e23101902}=&-\Big\{(A+B\hat{K}^N_1 -B\hat{K}^N_2 )^T[(\hat{P}^N_4 -\hat{P}^N_3 )\hat{x}_i +N\hat{P}^N_3 \hat{x}^{(N)} +\hat{s}^N_2 ]+(G+NB\hat{K}^N_2)^T\Big[[\frac{1}{N}(\hat{P}^N_1 -\hat{P}^N_2)\\&+\frac{N-1}{N}(\hat{P}^N_4 -\hat{P}_3) ]\hat{x}_i +[\hat{P}^N_2 +(N-1)\hat{P}^N_3]\hat{x}^{(N)} +\frac{1}{N}\hat{s}^N_1+\frac{N-1}{N}\hat{s}^N_2 \Big]+ [(\hat{K}^N_2)^TR(\hat{K}^N_1-\hat{K}^N_2)  -\frac{\Gamma^TQ}{N}]\hat{x}_i \nonumber\\&+[N(\hat{K}^N_2)^TR\hat{K}^N_2 +\frac{\Gamma^TQ\Gamma}{N}]\hat{x}^{(N)} +(\hat{K}^N_2)^TR \hat{v}^N    +\frac{\Gamma^TQ}{N}\eta\Big\}dt + (\hat{P}^N_4 -\hat{P}^N_3 )\sigma dw_i +\sum_{i=1}^N\hat{P}^N_3 \sigma dw_i .\nonumber
						\end{align}
					\end{subequations}
					By the uniqueness of the solutions of \eqref{e23101703} and \eqref{e23110701}, we have the following correspondence for $t\in[0,T], \forall i\in\mathcal{N},\forall j,k\in\mathcal{N}_{-i},$
					\begin{subequations}\label{e101605}
						\begin{align}
							\hat{p}^i_i =&\ (\hat{P}^N_1 -\hat{P}^N_2)\hat{x}_i +N\hat{P}^N_2 \hat{x}^{(N)} +\hat{s}^N_1 ,\\
							\hat{p}^i_j =&\ (\hat{P}^N_4 -\hat{P}^N_3 )\hat{x}_i +N\hat{P}^N_3 \hat{x}^{(N)} +\hat{s}^N_2 ,\\
							\begin{split}
								\hat{p}^{i,(N)} =&\ \Big[\frac{1}{N}(\hat{P}^N_1 -\hat{P}^N_2 )+\frac{N-1}{N}(\hat{P}^N_4 -\hat{P}_3 )\Big]\hat{x}_i +[\hat{P}^N_2 +(N-1)\hat{P}^N_3 ]\hat{x}^{(N)} +\frac{1}{N}\hat{s}^N_1 +\frac{N-1}{N}\hat{s}^N_2 ,
							\end{split} \\
							\hat{z}^i_{ii} =&\ \hat{P}^N_1 \sigma,\quad \hat{z}^i_{ij} =\hat{P}^N_2 \sigma,\quad \hat{z}^i_{jk} =\hat{P}^N_3 \sigma,\quad 
							\hat{z}^i_{ji} =\hat{P}^N_4 \sigma.
						\end{align} 
					\end{subequations}
					Substituting \eqref{e101605} into the stationarity condition \eqref{e23101801}, we obtain  
					\begin{equation} 
						\begin{split}
							0=&\	(B^T\hat{P}^N_1 +R\hat{K}^N_1)\hat{x}_i +\sum_{j\neq i}(B^T\hat{P}^N_2 +R\hat{K}^N_2 )\hat{x}_j +B^T\hat{s}^N_1 +R\hat{v}^N .
						\end{split}
					\end{equation}
					The representation \eqref{e101701} is obtained by the arbitrariness of $\hat{x}_i,i\in\mathcal{N}$. Furthermore, it follows by substituting
					\eqref{e101701} into \eqref{e101601} and \eqref{e101602} that
					\begin{subequations} 
						\begin{align}
							\begin{split} 
								\left\{
								\begin{aligned}
									&\dot{\hat{P}}^N_1 = \hat{P}^N_1 \Upsilon \hat{P}^N_1 +(N-1)\hat{P}^N_2 \Upsilon \hat{P}^N_2 +(N-1)(\hat{P}^N_2)^T\Upsilon \hat{P}^N_4-\hat{P}^N_1 (A+\frac{G}{N})\\&  \quad\quad\quad-(A+\frac{G}{N})^T\hat{P}^N_1 \!-\!\frac{N\!-\!1}{N}\hat{P}^N_2 G\!-\!\frac{N\!-\!1}{N}G^T\hat{P}^N_4 \!-\!Q\!+\!\frac{Q\Gamma}{N}\!+\!\frac{\Gamma^TQ}{N}\!-\!\frac{\Gamma^TQ\Gamma}{N^2}, \\&\hat{P}^N_1(T)=Q_f-\frac{ \Gamma_f^TQ_f}{N}-\frac{Q_f\Gamma_f}{N}+\frac{\Gamma_f^TQ_f\Gamma_f}{N^2},
								\end{aligned}
								\right.
							\end{split}\\
							\begin{split} 
								\left\{
								\begin{aligned}
									&\dot{\hat{P}}^N_2 =\hat{P}^N_1 \Upsilon \hat{P}^N_2 +(N-1)(\hat{P}^N_2)^T\Upsilon \hat{P}^N_3 +\hat{P}^N_2 \Upsilon \hat{P}^N_1 +(N-2)\hat{P}^N_2 \Upsilon \hat{P}^N_2-\hat{P}^N_1 \frac{G}{N}\\& \quad\quad\quad -\hat{P}^N_2 (A+\frac{N\!-\!1}{N}G)-(A+\frac{G}{N})^T\hat{P}^N_2 -\frac{N\!-\!1}{N}G^T\hat{P}^N_3 -\frac{\Gamma^TQ\Gamma}{N^2}+\frac{Q\Gamma}{N} , 	\\&\hat{P}^N_2(T)=\frac{\Gamma_f^TQ_f\Gamma_f}{N^2}-\frac{Q_f\Gamma_f}{N},
								\end{aligned}
								\right.
							\end{split}\\
							\begin{split} 
								\left\{
								\begin{aligned}
									&\dot{\hat{P}}^N_3 =(\hat{P}^N_1)^T\Upsilon \hat{P}^N_3 \!+\!\hat{P}^N_3 \Upsilon \hat{P}^N_1  \!+\!(N\!-\!2)\hat{P}^N_3 \Upsilon \hat{P}^N_2 \!+\!\hat{P}^N_4 \Upsilon \hat{P}^N_2 \!+\!(N\!-\!2)(\hat{P}^N_2)^T\Upsilon \hat{P}^N_3\!\!\\&\quad\quad\quad -\hat{P}^N_4 \frac{G}{N}-\frac{G^T}{N}\hat{P}^N_2  -\hat{P}^N_3 (A+\frac{N-1}{N}G) -(A+\frac{N-1}{N}G)^T\hat{P}^N_3  -\frac{\Gamma^TQ\Gamma}{N^2}, \!\!\!\!\!\! \\	&\hat{P}^N_3(T)=\frac{\Gamma_f^TQ_f\Gamma_f}{N^2},
								\end{aligned}
								\right.
							\end{split}\\
							\begin{split} 
								\left\{
								\begin{aligned}
									&\dot{\hat{P}}^N_4 =\hat{P}^N_4 \Upsilon \hat{P}^N_1  \!+\!(\hat{P}^N_1)^T\Upsilon \hat{P}^N_4\!+\!(N\!-\!1)\hat{P}^N_3 \Upsilon \hat{P}^N_2 \!+\!(N\!-\!2)(\hat{P}^N_2)^T\Upsilon \hat{P}^N_4\! -\!\frac{G^T}{N}\hat{P}^N_1 \!\! \\& \quad\quad\quad -	\hat{P}^N_4 (A+\frac{G}{N})-\frac{N\!-\!1}{N}\hat{P}^N_3 G -(A+\frac{N\!-\!1}{N}G)^T\hat{P}^N_4 +\frac{\Gamma^TQ}{N} -\frac{\Gamma^TQ\Gamma}{N^2}, 	\\&\hat{P}^N_4(T)=\frac{\Gamma_f^TQ_f\Gamma_f}{N^2}-\frac{\Gamma_f^TQ_f}{N},
								\end{aligned}
								\right.
							\end{split}\\
							\begin{split}
								\left\{
								\begin{aligned}
									&\dot{\hat{s}}^N_1 =-[A+\frac{G}{N}-\Upsilon (\hat{P}^N_1)^T-(N-1)\Upsilon (\hat{P}^N_2)^T]^T\hat{s}^N_1  -(N-1)(\frac{G}{N}-\Upsilon \hat{P}^N_2)^T\hat{s}^N_2 +(Q-\frac{\Gamma^TQ}{N})\eta, \\&\hat{s}^N_1(T)=(\frac{\Gamma_f^TQ_f}{N}-Q_f)\eta_f,
								\end{aligned}
								\right.
							\end{split}\\
							\begin{split}
								\left\{
								\begin{aligned}
									& \dot{\hat{s}}^N_2 =-[\frac{G}{N} -\Upsilon (\hat{P}^N_4)^T-(N-1)\Upsilon (\hat{P}^N_3)^T]^T\hat{s}^N_1  -\frac{\Gamma^TQ}{N}\eta-[A+\frac{N-1}{N}G-\Upsilon \hat{P}^N_1 -(N-2)\Upsilon \hat{P}^N_2 ]^T\hat{s}^N_2  ,\\	&\hat{s}^N_2(T)=\frac{\Gamma_f^TQ_f}{N}\eta_f .
								\end{aligned}
								\right.
							\end{split}
						\end{align}	
					\end{subequations}
					By computation, we have $(\hat{P}^N_1 )^T=\hat{P}^N_1 $, $(\hat{P}^N_2) ^T=\hat{P}^N_4 $, $(\hat{P}^N_3 )^T=\hat{P}^N_3 $, which further gives the Riccati equations \eqref{e101409} and ODEs \eqref{e101410}. This completes the proof of necessity.

					The sufficiency is proved below. We take any $u_i(\cdot)\in\mathcal{U}_c[0,T], \forall i\in\mathcal{N}$. By \eqref{e100101}, the dynamics of $x^{(N)}(\cdot)$ is given by
					\begin{equation}\label{e101702}
						\begin{split}
							dx^{(N)} =[(A+G)x^{(N)} +Bu^{(N)} ]dt+\sum_{i=1}^N\frac{\sigma}{N}dw_i ,\quad x^{(N)}(0)=x^{(N)}_0,
						\end{split}
					\end{equation}
					where $u^{(N)}(\cdot)\triangleq (1/N)\sum_{i=1}^Nu_i(\cdot)$. By \eqref{e100101}, \eqref{e101702}, \eqref{e101409}, \eqref{e101410}, and applying It\^{o}'s formula to $	t\mapsto \mathbb{E}	\langle x_i , [\hat{P}^N_1 -\hat{P}^N_2 -(\hat{P}^N_2)^T+\hat{P}^N_3 ] x_i \rangle$, $t\mapsto	\mathbb{E}	\langle x_i ,N(\hat{P}^N_2 -\hat{P}^N_3 )x^{(N)} \rangle$, $t\mapsto	\mathbb{E} \langle x^{(N)} ,N^2\hat{P}^N_3 x^{(N)} \rangle$, $t\mapsto	\mathbb{E}	\langle \hat{s}^N_1 -\hat{s}^N_2 ,x_i \rangle$ and $t\mapsto	\mathbb{E}	\langle N\hat{s}^N_2 ,x^{(N)} \rangle$,
					we have
					\begin{subequations}\label{e101703}
						\begin{align}
							\begin{split}
								&\ \mathbb{E}\langle x_i(T), Q_f x_i(T)\rangle-\mathbb{E}\langle x_{i0}, [\hat{P}^N_1(0)-\hat{P}^N_2(0)-\hat{P}^N_2(0)^T+\hat{P}^N_3(0)] x_{i0}\rangle\\
								=&\ \mathbb{E}\int_0^T\Big\{2\langle Gx^{(N)} +Bu_i ,[\hat{P}^N_1 -\hat{P}^N_2 -(\hat{P}^N_2)^T+\hat{P}^N_3 ]x_i \rangle-\langle x_i ,\big\{-[\hat{P}^N_1 -(\hat{P}^N_2)^T]\Upsilon (\hat{P}^N_1 -\hat{P}^N_2 )\\&-(-\hat{P}^N_2 +\hat{P}^N_3 )\Upsilon (\hat{P}^N_1 -\hat{P}^N_2 )-[\hat{P}^N_1 -(\hat{P}^N_2)^T]\Upsilon [-(\hat{P}^N_2)^T+\hat{P}^N_3 ]+Q\big\}x_i \rangle\\&+\langle \sigma, [\hat{P}^N_1 -\hat{P}^N_2 -(\hat{P}^N_2)^T+\hat{P}^N_3 ]\sigma\rangle \Big\}dt,
							\end{split}\\
							\begin{split}
								&\ \mathbb{E}\langle x_i(T),-Q_f\Gamma_fx^{(N)}(T)\rangle-\mathbb{E}\langle x_{i0},N[\hat{P}^N_2(0)-\hat{P}^N_3(0)]x^{(N)}_0\rangle\\
								=&\ \mathbb{E}\int_0^T\Big\{\langle Gx^{(N)} ,N(\hat{P}^N_2 -\hat{P}^N_3 )x^{(N)} \rangle+\langle Bu_i ,N[\hat{P}^N_2 -\hat{P}^N_3 ]x^{(N)} \rangle-\langle x_i ,\big\{[\hat{P}^N_1 -\hat{P}^N_2 -(\hat{P}^N_2)^T+\hat{P}^N_3 ]G\\&-(\hat{P}^N_2 -\hat{P}^N_3)\Upsilon N(\hat{P}^N_1 -\hat{P}^N_2)-[\hat{P}^N_1 \!-\!(\hat{P}^N_2)^T]\Upsilon N(\hat{P}^N_2 \!-\!\hat{P}^N_3 )-(N\!-\!1)(\hat{P}^N_2 -\hat{P}^N_3 )\Upsilon N\hat{P}^N_2 -Q\Gamma\big\}x^{(N)} \rangle \\&+\langle x_i ,N(\hat{P}^N_2 -\hat{P}^N_3 )Bu^{(N)} \rangle+\langle \sigma,(\hat{P}^N_2 -\hat{P}^N_3 )\sigma\rangle \Big\}dt, \!\!\!\!\!
							\end{split}\\
							\begin{split}
								&\ \mathbb{E}\langle x^{(N)}(T),\Gamma^T_fQ_f\Gamma_fx^{(N)}(T)\rangle -\mathbb{E}\langle x^{(N)}_0,N^2\hat{P}^N_3(0)x^{(N)}_0\rangle\\
								=&\ \mathbb{E}\int_0^T\Big\{\langle Bu^{(N)} , N^2\hat{P}^N_3 x^{(N)}  \rangle-\langle x^{(N)} , \big\{N[(\hat{P}^N_2)^T-\hat{P}^N_3 ]G-N^2[(\hat{P}^N_2)^T-\hat{P}^N_3 ]\Upsilon \hat{P}^N_2 -N^2\hat{P}^N_3 \Upsilon \hat{P}^N_1 \\&-N^2\hat{P}^N_3 (N-1)\Upsilon \hat{P}^N_2 -[\hat{P}^N_1 -(\hat{P}^N_2)^T]\Upsilon N^2\hat{P}^N_3 +G^T(N\hat{P}^N_2 -N\hat{P}^N_3)-N(\hat{P}^N_2)^T\Upsilon  N(N-1)\hat{P}^N_3 \\& +\Gamma^TQ\Gamma\big\}x^{(N)}  \rangle +\langle x^{(N)} , N^2\hat{P}^N_3 Bu^{(N)}  \rangle+\langle \sigma ,  N\hat{P}^N_3 \sigma \rangle \Big\}dt,
							\end{split}\\
							\begin{split}
								&\ \mathbb{E}\langle -Q_f\eta_f,x_i(T)\rangle-\mathbb{E}\langle \hat{s}^N_1(0)-\hat{s}^N_2(0),x_{i0}\rangle	\\
								=&\ \mathbb{E}\int_0^T\Big\{-\langle -[\hat{P}^N_1 -\hat{P}^N_2 -(\hat{P}^N_2)^T+\hat{P}^N_3 ]\Upsilon \hat{s}^N_1 -N(\hat{P}^N_2 -\hat{P}^N_3)\Upsilon \hat{s}^N_1 +[\hat{P}^N_1 -(\hat{P}^N_2)^T]\Upsilon \hat{s}^N_2 -Q\eta,x_i \rangle \\&+\langle \hat{s}^N_1 -\hat{s}^N_2 ,Gx^{(N)} \rangle+\langle \hat{s}^N_1 \!-\!\hat{s}^N_2 ,Bu_i  \rangle\Big\}dt,
							\end{split}\\
							\begin{split}
								&\ \mathbb{E}\langle \Gamma_f^TQ_f\eta_f ,x^{(N)}(T)\rangle-\mathbb{E}\langle N\hat{s}^N_2(0),x^{(N)}_0\rangle\\
								=&\ \mathbb{E}\int_0^T\!\Big\{\!-\langle N\hat{P}^N_3 \Upsilon \hat{s}^N_1\! -\!N^2\hat{P}^N_3 \Upsilon \hat{s}^N_1 \!-\![\hat{P}^N_1\! -\!(\hat{P}^N_2)^T]\Upsilon N\hat{s}^N_2 +G^T(\hat{s}^N_1 -\hat{s}^N_2 )-N(\hat{P}^N_2)^T\Upsilon [(N-1)\hat{s}^N_2 \\&+\hat{s}^N_1  ]+\Gamma^TQ\eta,x^{(N)} \rangle +\langle N\hat{s}^N_2 ,Bu^{(N)} \rangle \Big\}dt.
							\end{split}
						\end{align}
					\end{subequations}
					Substituting the terminal terms  in the cost functional \eqref{e100102} by using \eqref{e101703}, we obtain 
					\begin{equation}
						\begin{split}		
							&\	{\rm J}_i(\mathbf{x}_0;u_i(\cdot),u_j(\cdot),\forall j\in\mathcal{N}_{-i})\\
							=&\ |\eta_f|^2_{Q_f}+\mathbb{E}[\langle x_{i0}, [\hat{P}^N_1(0)-2\hat{P}^N_2(0)+\hat{P}^N_3(0)] x_{i0}\rangle+2\langle  x_{i0},N[\hat{P}^N_2(0)-\hat{P}^N_3(0)]x^{(N)}_0\rangle\\&+\langle  x^{(N)}_0,N^2\hat{P}^N_3(0) x^{(N)}_0\rangle+2\langle \hat{s}^N_1(0)-\hat{s}^N_2(0),x_{i0}\rangle+2\langle N\hat{s}^N_2(0),x^{(N)}_0\rangle]\\&+\mathbb{E}\int_0^T\Big[|u_i +R^{-1}B^T(\hat{P}^N_1 -\hat{P}^N_2)x_i +NR^{-1}B^T\hat{P}^N_2 x^{(N)} +R^{-1}B^T\hat{s}^N_1 |^2_{R}
							\\&+\langle \sigma, [\hat{P}^N_1 +(N-1)\hat{P}^N_3 ]\sigma\rangle   -\langle 2(N-1)\Upsilon \hat{s}^N_2 + \Upsilon \hat{s}^N_1 ,\hat{s}^N_1 \rangle	+	|\eta|^2_Q \Big]dt
						\end{split}
					\end{equation}
					\begin{equation}
						\begin{split}		
							{\rm J}_i(\mathbf{x}_0;u_i(\cdot),u_j(\cdot),j\in\mathcal{N}_{-i})
							= V_i(\mathbf{x}_0)
							+\mathbb{E}\int_0^T\big|u_i -(\hat{K}^N_1 -\hat{K}^N_2)x_i -N\hat{K}^N_2 x^{(N)} -\hat{v}^N \big|^2_{R}dt,
						\end{split}
					\end{equation}
					where the representation \eqref{e101701} is used. Thus \eqref{e101402} holds if and only if  
					$$\hat{u}_i =(\hat{K}^N_1 -\hat{K}^N_2 )\hat{x}_i +N\hat{K}^N_2 \hat{x}^{(N)} +\hat{v}^N.
					$$
					This shows that $(\hat{K}^N_1(\cdot),\hat{K}^N_2(\cdot),\hat{v}^N(\cdot))$ is a centralized closed-loop Nash equilibrium of Problem (P2). The proof is completed.
					
					\section{Proofs for Results in Section \ref{decentralized}}
					\subsection{Proof of Theorem \ref{theorem101901}}\label{proofoftheorem101901}
					It follows by \eqref{e101201} and \eqref{e042401} that   
					\begin{subequations}\label{e101801}
						\begin{align}
							\begin{split}
								\left\{
								\begin{aligned}
									&\dot{\check{\Lambda}}^N_1 =\check{\Lambda}^N_1 \Upsilon \check{\Lambda}^N_1 -\check{\Lambda}^N_1 A-A^T\check{\Lambda}^N_1 -Q+\check{g}_1(t,\frac{1}{N}),\\& 	\check{\Lambda}^N_1(T)=Q_f+\check{g}_1(T,\frac{1}{N}),
								\end{aligned}
								\right.
							\end{split}\\
							\begin{split}
								\left\{
								\begin{aligned}
									&\dot{\check{\Lambda}}^N_2 =\check{\Lambda}^N_1 \Upsilon \check{\Lambda}^N_2 \!+\!\check{\Lambda}^N_2 \Upsilon \check{\Lambda}^N_1 \!+\!\check{\Lambda}^N_2 \Upsilon \check{\Lambda}^N_2 \!-\!\check{\Lambda}^N_1 G\!-\!\check{\Lambda}^N_2 (A\!+\!G)\!-\!A^T\check{\Lambda}^N_2 \!+\!Q\Gamma\!+\!\check{g}_2(t,\frac{1}{N}),\\& 	\check{\Lambda}^N_2(T)=-Q_f\Gamma_f+\check{g}_2(T,\frac{1}{N}),
								\end{aligned}
								\right.
							\end{split}\\
							\begin{split}
								\left\{
								\begin{aligned}
									&\dot{\check{\Lambda}}^N_3 =\check{\Lambda}^N_4 \Upsilon \check{\Lambda}^N_2 \!+\!\check{\Lambda}^N_3 \Upsilon \check{\Lambda}^N_1 \!+\!\check{\Lambda}^N_3 \Upsilon \check{\Lambda}^N_2 \!-\!\check{\Lambda}^N_3 (A\!+\!G)\!-\!\check{\Lambda}^N_4 G\!-\!(A\!+\!G)^T\check{\Lambda}^N_3 \!-\!G^T\check{\Lambda}^N_2 \!-\!\Gamma^TQ\Gamma\!+\!\check{g}_3(t,\frac{1}{N}),\quad 	\\&\check{\Lambda}^N_3(T)=\Gamma^T_fQ_f\Gamma_f+\check{g}_3(T,\frac{1}{N}),
								\end{aligned}
								\right.
							\end{split}\\
							\begin{split}
								\left\{
								\begin{aligned}
									&\dot{\check{\Lambda}}^N_4 =\check{\Lambda}^N_4 \Upsilon \check{\Lambda}^N_1  -\check{\Lambda}^N_4 A-(A+G)^T\check{\Lambda}^N_4 -G^T\check{\Lambda}^N_1 +\Gamma^TQ+\check{g}_4(t,\frac{1}{N}),\quad 	\\&\check{\Lambda}^N_4(T)=-\Gamma^T_fQ_f+\check{g}_4(T,\frac{1}{N}),
								\end{aligned}
								\right.
							\end{split}\\
							\begin{split}
								\left\{
								\begin{aligned}
									&\dot{\check{\varphi}}^N_1 =-(A^T-\check{\Lambda}^N_1 \Upsilon -\check{\Lambda}^N_2 \Upsilon )\check{\varphi}^N_1 +Q\eta+\check{g}_5(t,\frac{1}{N}),\quad \\&\check{\varphi}^N_1(T)=-Q_f\eta_f+\check{g}_5(T,\frac{1}{N}),
								\end{aligned}
								\right.
							\end{split}\\
							\begin{split}
								\left\{
								\begin{aligned}
									&\dot{\check{\varphi}}^N_2 =-(G^T-\check{\Lambda}^N_4 \Upsilon +\check{\Lambda}^N_3 \Upsilon )\check{\varphi}^N_1 -(A+G)^T\check{\varphi}^N_2 \Gamma^TQ\eta+\check{g}_6(T,\frac{1}{N}),\quad\\ 	&\check{\varphi}^N_2(T)=\Gamma_f^TQ_f\eta_f+\check{g}_6(T,\frac{1}{N}), 
								\end{aligned}
								\right.
							\end{split}
						\end{align}	
					\end{subequations}
					where
					\begin{subequations}\label{e013101}
						\begin{align}
							\begin{split}
								\left\{
								\begin{aligned}
									&\check{g}_1(t,\frac{1}{N})=\frac{N-1}{N^2}\check{\Lambda}^N_2 \Upsilon \check{\Lambda}^N_2 -\check{\Lambda}^N_1 \frac{G}{N}-\frac{G^T}{N}\check{\Lambda}^N_1 -\frac{N-1}{N^2}\check{\Lambda}^N_2 G-\frac{N-1}{N^2}G^T\check{\Lambda}^N_4 +\frac{\Gamma^TQ}{N}+\frac{Q\Gamma}{N}-\frac{\Gamma^TQ\Gamma}{N^2},\\
									&\check{g}_1(T,\frac{1}{N})=-\frac{\Gamma^T_fQ_f}{N}-\frac{Q_f\Gamma_f}{N}+\frac{\Gamma^T_fQ_f\Gamma_f}{N^2},
								\end{aligned}
								\right.
							\end{split}\\
							\begin{split}
								\left\{
								\begin{aligned}
									&\check{g}_2(t,\frac{1}{N})=-\frac{2}{N}\check{\Lambda}^N_2 \Upsilon \check{\Lambda}^N_2 +\check{\Lambda}^N_2 \frac{G}{N}-\frac{G^T}{N}\check{\Lambda}^N_2-\frac{N-1}{N^2}G^T\check{\Lambda}^N_3 -\frac{\Gamma^TQ\Gamma}{N},\\& 	\check{g}_2(T,\frac{1}{N})=\frac{\Gamma^T_fQ_f\Gamma_f}{N},
								\end{aligned}
								\right.
							\end{split}\\
							\begin{split}
								\left\{
								\begin{aligned}
									&\check{g}_3(t,\frac{1}{N})=-\frac{2}{N}\check{\Lambda}^N_3 \Upsilon \check{\Lambda}^N_2 +\check{\Lambda}^N_3 \frac{G}{N}-\frac{G^T}{N}\check{\Lambda}^N_3 ,\quad 	\\&\check{g}_3(T,\frac{1}{N})=0,
								\end{aligned}
								\right.
							\end{split}\\
							\begin{split}
								\left\{
								\begin{aligned}
									&\check{g}_4(t,\frac{1}{N})=\frac{N-1}{N^2}\check{\Lambda}^N_3 \Upsilon \check{\Lambda}_2 -\check{\Lambda}^N_4 \frac{G}{N}+\frac{G^T}{N}	\check{\Lambda}^N_4 -\frac{N-1}{N^2}\check{\Lambda}^N_3 G-\frac{\Gamma^TQ\Gamma}{N},\quad 	\\&\check{g}_4(T,\frac{1}{N})=\frac{\Gamma^T_fQ_f\Gamma_f}{N},
								\end{aligned}
								\right.
							\end{split}\\
							\begin{split}
								\left\{
								\begin{aligned}
									&\check{g}_5(t,\frac{1}{N})=-\frac{1}{N}\check{\Lambda}^N_2 \Upsilon \check{\varphi}^N_1 -\frac{G^T}{N}\check{\varphi}^N_1  -\frac{N-1}{N^2}G^T\check{\varphi}^N_2 -\frac{\Gamma^TQ}{N}\eta,\quad \\&\check{g}_5(T,\frac{1}{N})=\frac{\Gamma_f^TQ_f}{N}\eta_f,
								\end{aligned}
								\right.
							\end{split}\\
							\begin{split}
								\left\{
								\begin{aligned}
									&\check{g}_6(t,\frac{1}{N})=-\frac{1}{N}\check{\Lambda}^N_3 \Upsilon \check{\varphi}^N_1 +\frac{G^T}{N}\check{\varphi}^N_2 ,\quad\\ 	&\check{g}_6(t,\frac{1}{N})=0.
								\end{aligned}
								\right.
							\end{split}
						\end{align}	
					\end{subequations}
					The proof of necessity is completed by	letting $N\rightarrow\infty$ in \eqref{e101801}. In view of $\check{g}_k,k=1,2,. . .,6$, the sufficiency follows from \cite[Theorem 4, i)]{b110205}.
					
					\subsection{Proof of Theorem \ref{theorem101902}}\label{proofoftheorem101902}
					Similar to \eqref{e101801}, it follows by	\eqref{e101409}, \eqref{e101410} that
					\begin{subequations}\label{e101802}
						\begin{align}
							\begin{split} 
								\left\{
								\begin{aligned}
									&\dot{\hat{\Lambda}}^N_1 =\hat{\Lambda}^N_1 \Upsilon\hat{\Lambda}^N_1 -\hat{\Lambda}^N_1 A-A^T\hat{\Lambda}^N_1 -Q+\hat{g}_1(t,\frac{1}{N}), \\&\hat{\Lambda}^N_1(T)=Q_f+\hat{g}_1(T,\frac{1}{N}),
								\end{aligned}
								\right.
							\end{split}\\
							\begin{split} 
								\left\{
								\begin{aligned}
									&\dot{\hat{\Lambda}}^N_2 =\hat{\Lambda}^N_1 \Upsilon\hat{\Lambda}^N_2 +\hat{\Lambda}^N_2 \Upsilon\hat{\Lambda}^N_1 +\hat{\Lambda}^N_2 \Upsilon\hat{\Lambda}^N_2 -\hat{\Lambda}^N_1 G-\hat{\Lambda}^N_2 (A+G) -A^T\hat{\Lambda}^N_2 +Q\Gamma+\hat{g}_2(t,\frac{1}{N}) ,\quad 	\\&\hat{\Lambda}^N_2(T)=-Q_f\Gamma_f+\hat{g}_2(T,\frac{1}{N}),
								\end{aligned}
								\right.
							\end{split}\\
							\begin{split} 
								\left\{
								\begin{aligned}
									&\dot{\hat{\Lambda}}^N_3 =\hat{\Lambda}^N_1 \Upsilon\hat{\Lambda}^N_3 +\hat{\Lambda}^N_3 \Upsilon\hat{\Lambda}^N_1 +(\hat{\Lambda}^N_2)^T\Upsilon\hat{\Lambda}^N_2 +\hat{\Lambda}^N_3 \Upsilon\hat{\Lambda}^N_2 +(\hat{\Lambda}^N_2)^T\Upsilon\hat{\Lambda}^N_3 \\&\quad\quad-\!(\hat{\Lambda}^N_2)^TG\!-\!G^T\hat{\Lambda}^N_2 \!-\!\hat{\Lambda}^N_3 (A\!+\!G) \!-\!(A\!+\!G)^T\hat{\Lambda}^N_3  \!-\!\Gamma^TQ\Gamma\!+\!\hat{g}_3(t,\frac{1}{N}),\quad \\	&\hat{\Lambda}^N_3(T)=\Gamma_f^TQ_f\Gamma_f+\hat{g}_3(T,\frac{1}{N}),
								\end{aligned}
								\right.
							\end{split}\\
							\begin{split}  
								\left\{
								\begin{aligned}
									&\dot{\hat{\varphi}}^N_1 =-[A-\Upsilon\hat{\Lambda}^N_1 -\Upsilon(\hat{\Lambda}^N_2) ^T]^T\hat{\varphi}^N_1 +Q\eta+\hat{g}_4(t,\frac{1}{N}), \\&\hat{\varphi}^N_1(T)=-Q_f\eta_f+\hat{g}_4(T,\frac{1}{N}),
								\end{aligned}
								\right.
							\end{split}\\
							\begin{split}  
								\left\{
								\begin{aligned}
									& \dot{\hat{\varphi}}^N_2 =-(G-\Upsilon\hat{\Lambda}^N_2 -\Upsilon\hat{\Lambda}^N_3 )^T\hat{\varphi}^N_1 -(A+G-\Upsilon\hat{\Lambda}^N_1 -\Upsilon\hat{\Lambda}^N_2)^T\hat{\varphi}^N_2  -\Gamma^TQ\eta+\hat{g}_5(t,\frac{1}{N}),\\	&\hat{\varphi}^N_2(T)=\Gamma_f^TQ_f\eta_f+\hat{g}_5(T,\frac{1}{N}).
								\end{aligned}
								\right.
							\end{split}
						\end{align}	
					\end{subequations}
					In particular, we can determine
					\begin{align}
						\begin{split} 
							\left\{
							\begin{aligned}
								&\hat{g}_1(t,\frac{1}{N})=\frac{N-1}{N^2}\hat{\Lambda}^N_2 \Upsilon\hat{\Lambda}^N_2 +\frac{N-1}{N^2}(\hat{\Lambda}^N_2)^T\Upsilon(\hat{\Lambda}^N_2)^T-\hat{\Lambda}^N_1 \frac{G}{N}+\frac{G^T}{N}\hat{\Lambda}^N_1\\&\quad\quad\quad\quad\quad -\frac{N-1}{N^2}\hat{\Lambda}^N_2 G-\frac{N-1}{N^2}G^T(\hat{\Lambda}^N_2)^T+\frac{Q\Gamma}{N}+\frac{\Gamma^TQ}{N}-\frac{\Gamma^TQ\Gamma}{N^2}, \\&\hat{g}_1(T,\frac{1}{N})=-\frac{ \Gamma_f^TQ_f}{N}-\frac{Q_f\Gamma_f}{N}+\frac{\Gamma_f^TQ_f\Gamma_f}{N^2},
							\end{aligned}
							\right.
						\end{split}
					\end{align}	 
					The expressions of $\hat{g}_2,\hat{g}_3,\hat{g}_4$ and $\hat{g}_5$ can be e determined in a similar way and the detail is omitted here. The proof of necessity is completed by letting $N\rightarrow\infty$ in \eqref{e101801}. In view of $\check{g}_k,k=1,2,. . .,6$, the sufficiency follows from \cite[Theorem 4, i)]{b110205}.
					
					
					\subsection{Proof of Proposition \ref{t122301}}\label{proofoft122301}
					
					By \eqref{e101903} and \eqref{e101902}, the dynamics of $ \check{x}^{(N)}-\bar{x}$ is given by
					\begin{equation} 
						\begin{split}
							d(\check{x}^{(N)}-\bar{x}) =&\ \big[(A+G-\Upsilon\check{\Lambda}_1^N -\Upsilon\check{\Lambda}_2^N)(\check{x}^{(N)}-\bar{x}) +\frac{1}{N}\Upsilon \check{\Lambda}_2^N \check{x}^{(N)}+\Upsilon(P_1^\infty-\check{\Lambda}_1^N +P_2^\infty-\check{\Lambda}_2^N) \bar{x}  \\&-\Upsilon (\check{\varphi}_1^N- s_1^\infty) \big]dt +\sum_{i=1}^N\frac{\sigma}{N} dw_i ,\quad \check{x}^{(N)}(0)-\bar{x}(0)=x^{(N)}_{0}-\bar{x}_0.
						\end{split}
					\end{equation}
					Similarly  we can obtain the dynamics of $\hat{x}^{(N)}-\bar{x}$ by using \eqref{e101906} and \eqref{e101902}. The proposition follows from elementary estimates by use of Theorem \ref{theorem101901} and Theorem \ref{theorem101902}.

\end{document}